\documentclass[11pt]{amsart}
\usepackage{palatino, mathpazo}
\usepackage[mathscr]{eucal}
\usepackage{amssymb}
\usepackage[all]{xy}
\usepackage{graphicx}
\usepackage{texdraw}

\newtheorem{theorem}{Theorem}[section]
\newtheorem{lemma}[theorem]{Lemma}

\newtheorem{corollary}[theorem]{Corollary}
\theoremstyle{remark}
\newtheorem{remark}[theorem]{Remark}
\newtheorem{example}[theorem]{Example}

\newcommand{\p}{\partial}
\newcommand{\T}{\vartheta}

\begin{document}

\title[Mean field equations on tori]
{Elliptic functions, Green functions and the mean field equations
on tori}
\date{June 1st, 2006. Revised June 29, 2008.}
\author{Chang-Shou Lin}
\address{Department of Mathematics and Taida Institute of
Mathematical Sciences (TIMS), Taipei.}
\email{cslin@math.ntu.edu.tw}
\author{Chin-Lung Wang}
\address{Department of Mathematics, National Taiwan University,
Taipei. Department of Mathematics, National Central University,
Jhongli.} \email{dragon@math.ntu.edu.tw; dragon@math.ncu.edu.tw}

\begin{abstract}
We show that the Green functions on flat tori can have either 3 or
5 critical points only. There does not seem to be any direct
method to attack this problem. Instead, we have to employ
sophisticated non-linear partial differential equations to study
it.

We also study the distribution of number of critical points over
the moduli space of flat tori through deformations. The functional
equations of special theta values provide important inequalities
which lead to a solution for all rhombus tori.
\end{abstract}

\maketitle

\section{Introduction and Statement of Results}

The study of geometric or analytic problems on two dimensional
tori is the same as the study of problems on $\mathbb{R}^2$ with
doubly periodic data. Such situations occur naturally in sciences
and mathematics since early days. The mathematical foundation of
elliptic functions was subsequently developed in the 19th century.
It turns out that these special functions are rather deep objects
by themselves. Tori of different shape may result in very
different behavior of the elliptic functions and their associated
objects. Arithmetic on elliptic curves is perhaps the eldest and
the most vivid example.

In this paper, we show that this is also the case for certain
non-linear partial differential equations. Indeed, researches on
doubly periodic problems in mathematical physics or differential
equations often restrict the study to rectangular tori for
simplicity. This leaves the impression that the theory for general
tori may resemble much the same way as for the rectangular case.
However this turns out to be false. We will show that the
solvability of the {\it mean field equation} depends on the shape
of the Green function, which in turn depends on the geometry of
the tori in an essential way.

Recall that the Green function $G(z, w)$ on a flat torus $T =
\mathbb{C}/\mathbb{Z}\omega_1 + \mathbb{Z}\omega_2$ is the unique
function on $T\times T$ which satisfies
\begin{equation*}
-\triangle_z G(z, w) = \delta_{w}(z) - \frac{1}{|T|}
\end{equation*}
and $\int_T G(z, w)\,dA = 0$, where $\delta_w$ is the Dirac
measure with singularity at $z = w$. Because of the translation
invariance of $\triangle_z$, we have $G(z,w) = G(z - w, 0)$ and it
is enough to consider {\it the Green function} $G(z) := G(z, 0)$.

Not surprisingly, $G$ can be explicitly solved in terms of
elliptic functions. For example, using theta functions we have
(cf.\ Lemma \ref{green-w}, Lemma \ref{green-theta})
$$
G(z) = -\frac{1}{2\pi}\log |\T_1(z)| + \frac{1}{2b}y^2 + C(\tau)
$$
where $z = x + iy$ and $\tau := \omega_2/\omega_1 = a + ib$. The
structure of $G$, especially its critical points and critical
values, will be the fundamental objects that interest us. The
critical point equation $\nabla G(z) = 0$ is given by
$$
\frac{\p G}{\p z} \equiv \frac{-1}{4\pi}\left((\log\T_1)_z + 2\pi
i \frac{y}{b}\right) = 0.
$$
In terms of Weierstrass' elliptic functions $\wp(z)$, $\zeta(z) :=
-\int^z \wp$ and using the relation $(\log\T_1)_z = \zeta(z) -
\eta_1 z$ with $\eta_i = \zeta(z + \omega_i) - \zeta(z)$ the
quasi-periods, the equation takes the simpler form: $z = t\omega_1
+ s\omega_2$ is a critical point of $G$ if and only if the
following linear relation (Lemma \ref{critical}) holds:
\begin{equation}\label{linear}
\zeta(t\omega_1 + s\omega_2) = t\eta_1 + s\eta_2.
\end{equation}

Since $G$ is even, it is elementary to see that half periods
$\frac{1}{2}\omega_1$, $\frac{1}{2}\omega_2$ and
$\frac{1}{2}\omega_3 = (\omega_1 + \omega_2)/2$ are the three
obvious critical points and other critical points must appear in
pair. The question is: {\it Are there other critical points?} or
{\it How many critical points might $G$ have?}. It turns out that
this is a delicate question and can not be attacked easily from
the simple looking equation (\ref{linear}). One of our chief
purposes in this paper is to understand the geometry of the
critical point set over the moduli space of flat tori
$\mathcal{M}_1 = \mathcal{H}/{\rm SL}(2,\mathbb{Z})$ and to study
its interaction with the non-linear mean field equation.\medskip

The mean field equation on a flat torus $T$ takes the form ($\rho
\in \Bbb R_+$)
\begin{equation}{\label{1-1}}
\triangle u + \rho e^u = \rho \delta_0.
\end{equation}
This equation has its origin in the prescribed curvature problem
in geometry like the Nirenberg problem, cone metrics etc.. It also
comes from statistical physics as the mean field limits of the
Euler flow, hence the name. Recently it was shown to be related to
the self dual condensation of Chern-Simons-Higgs model. We refer
to \cite{Chen-Lin1}, \cite{Chen-Lin2}, \cite{Chen-Lin3},
\cite{CLW}, \cite{JLW}, \cite{NT0}, \cite{NT} and \cite{NT1} for
the recent development of this subject.

When $\rho \ne 8m\pi$ for any $m\in\mathbb{Z}$, it has been
recently proved in \cite{Chen-Lin2}, \cite{Chen-Lin3}, \cite{CLW}
that the Leray-Schauder degree is non-zero, so the equation always
has solutions, regardless on the actual shape of $T$.

The first interesting case remained is when $\rho = 8\pi$ where
the degree theory fails completely. Instead of the topological
degree, the precise knowledge on the Green function plays a
fundamental role in the investigation of (\ref{1-1}). The first
main result of this paper is the following existence criterion
whose proof is given in \S3 by a detailed manipulation on elliptic
functions:

\begin{theorem}[Existence]{\label{main}}
For $\rho = 8\pi$, the mean field equation on a flat torus has
solutions if and only if the Green function has critical points
other than the three half period points. Moreover, each extra pair
of critical points corresponds to an one parameter scaling family
of solutions.
\end{theorem}

It is known that for rectangular tori $G(z)$ has precisely the
three obvious critical points, hence for $\rho = 8\pi$ equation
(\ref{1-1}) has no solutions. However we will show in \S2 that for
the case $\omega_1 = 1$ and $\tau = \omega_2 = e^{\pi i/3}$ there
are at least five critical points and the solutions of (\ref{1-1})
exist.

Our second main result is the uniqueness theorem.

\begin{theorem}[Uniqueness]\label{uniq}
For $\rho = 8\pi$, the mean field equation on a flat torus has at
most one solution up to scaling.
\end{theorem}

In view of the correspondence in Theorem \ref{main}, an equivalent
statement of Theorem \ref{uniq} is the following result:

\begin{theorem}\label{Green}
The Green function has at most five critical points.
\end{theorem}

Unfortunately we were unable to find a direct proof of Theorem
\ref{Green} from the critical point equation (\ref{linear}).
Instead, we will prove the uniqueness theorem first, and then
Theorem \ref{Green} is an immediate corollary. Our proof of
Theorem \ref{uniq} is based on the method of symmetrization
applied to the linearized equation at a particularly chosen even
solution in the scaling family. In fact we study in \S4 the one
parameter family
\begin{equation*}
\triangle u + \rho e^u = \rho \delta_0, \quad \rho \in [4\pi,
8\pi]
\end{equation*}
on $T$ within {\it even solutions}. This extra assumption allows
us to construct a double cover $T \to S^2$ via the Weierstrass
$\wp$ function and to transform equation (\ref{1-1}) into a
similar one on $S^2$ but with three more delta singularities with
negative coefficients. The condition $\rho \ge 4\pi$ is to
guarantee that the original singularity at $0$ still has
non-negative coefficient of delta singularity.

The uniqueness is proved for this family via the method of
continuity. For the starting point $\rho = 4\pi$, by a
construction similar to the proof of Theorem \ref{main} we sharpen
the result on nontrivial degree to the existence and uniqueness of
solution (Theorem \ref{4pi}). For $\rho \in [4\pi, 8\pi]$, the
symmetrization reduces the problem on the non-degeneracy of the
linearized equation to the isoperimetric inequality on domains in
$\mathbb{R}^2$ with respect to certain singular measure:

\begin{theorem}[Symmetrization Lemma]\label{sym}
Let $\Omega \subset \mathbb{R}^2$ be a simply connected domain and
let $v$ be a solution of
$$
\triangle v + e^v = \sum\nolimits_{j = 1}^N
2\pi\alpha_j\delta_{p_j}
$$
in $\Omega$. Suppose that the first eigenvalue of $\triangle +
e^v$ is zero on $\Omega$ with $\varphi$ the first eigenfunction.
If the isoperimetric inequality with respect to $ds^2 =
e^v|dx|^2$:
\begin{equation}\label{iso}
2 \ell^2 (\p\omega) \geq m(\omega) (4\pi - m(\omega))
\end{equation}
holds for all level domains $\omega = \{\varphi > t\}$ with $t >
0$, then
\begin{equation*}
\int_\Omega  e^v\, dx \geq 2\pi.
\end{equation*}
Moreover, (\ref{iso}) holds if there is only one negative
$\alpha_j$ and $\alpha_j = -1$.
\end{theorem}

The proof on the number of critical points appears to be one of
the very few instances that one needs to study a simple analytic
equation, here the critical point equation (\ref{linear}), by way
of sophisticated non-linear analysis.

To get a deeper understanding of the underlying structure of
solutions, we first notice that for $\rho = 8\pi$, (\ref{1-1}) is
the Euler-Lagrange equation of the non-linear functional
\begin{equation}\label{functional}
J_{8\pi}(v) = \frac{1}{2}\int_T |\nabla v|^2\,dA - 8\pi\log \int_T
e^{v - 8\pi G(z)}\,dA
\end{equation}
on $H^1(T) \cap \{v \mid \int_T v = 0 \}$, the Sobolev space of
functions with $L^2$-integrable first derivative. From this
viewpoint, the non-existence of minimizers for rectangular tori
was known in \cite{CLW}. Here we sharpen the result to the
non-existence of solutions. Also for $\rho \in (4\pi, 8\pi)$ we
sharpen the result on non-trivial degree of equation (\ref{1-1})
in \cite{Chen-Lin2} to the uniqueness of solutions within even
functions. We expect the uniqueness holds true without the even
assumption, but our method only achieves this at $\rho = 4\pi$.
Obviously, uniqueness without even assumption fails at $\rho =
8\pi$ due to the existence of scaling.

Naturally, the next question after Theorem \ref{uniq} is to
determine those tori whose Green functions have five critical
points. It is the case if the three half-periods are all saddle
points. A strong converse is proved in \cite{LW}:
\medskip

\noindent {\bf Theorem A.} {\it If the Green function has five
critical points then the extra pair of critical points are minimum
points and the three half-periods are all saddle points.}\medskip

Together with Theorem \ref{main}, this implies that {\it a
minimizer of $J_{8\pi}$ exists if and only if the Green function
has more than three critical points}. In fact we show in \cite{LW}
that any solution of equation (\ref{1-1}) must be a minimizer of
the nonlinear functional $J_{8\pi}$. Thus we completely solve the
existence problem on minimizers, a question raised by Nolasco and
Tarantello in \cite{NT}.

By Theorem A, we have reduced the question on detecting a given
torus to have five critical points to the technically much simpler
criterion on (non)-local minimality of the three half period
points. In this paper, however, no reference to Theorem A is
needed. Instead, it motivates the following comparison result,
which also simplifies the criterion further:

\begin{theorem}\label{comparison-intro}
Let $z_0$ and $z_1$ be two half period points. Then
$$
G(z_0) \ge G(z_1) \quad\mbox{if and only if}\quad |\wp(z_0)| \ge
|\wp(z_1)|.
$$
\end{theorem}

\begin{figure}
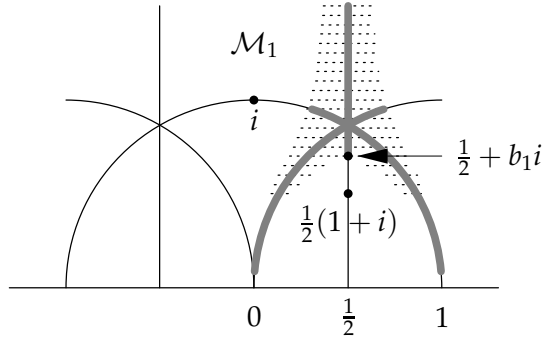

\renewcommand{\figurename}{Figure}
\begin{center}
\begin{texdraw}
  \drawdim cm  \setunitscale 2.5 \linewd 0.007
  \move(-1.3 0) \lvec(1.3 0)
  \move(0 0) \larc r:1 sd:0 ed:180
  \move(-0.5 0) \lvec (-0.5 1.5)
  \move(-1 0) \larc r:1 sd:0 ed:90
  \move(0.5 0) \lvec (0.5 1.5)
  \move(1 0) \larc r:1 sd:90 ed:180
  \textref h:C v:C \htext(0 1.3){$\mathcal{M}_1$}
  \textref h:C v:C \htext(0 -0.15){$0$}
  \textref h:C v:C \htext(0.5 -0.15){$\frac{1}{2}$}
  \textref h:C v:C \htext(1 -0.15){1}
  \move(0 1) \fcir f:0 r:0.025
  \textref h:C v:C \htext(0 0.9){$i$}
  \move(0.5 0.5) \fcir f:0 r:0.025
  \move(0.37 1.5) \lpatt(0.01 0.03) \lvec(0.63 1.5)
  \move(0.365 1.45) \lvec(0.635 1.45)
  \move(0.36 1.4) \lvec(0.64 1.4)
  \move(0.355 1.35) \lvec(0.645 1.35)
  \move(0.35 1.3) \lvec(0.65 1.3)
  \move(0.345 1.25) \lvec(0.655 1.25)
  \move(0.34 1.2) \lvec(0.66 1.2)
  \move(0.335 1.15) \lvec(0.665 1.15)
  \move(0.33 1.1) \lvec(0.67 1.1)
  \move(0.31 1.05) \lvec(0.69 1.05)
  \move(0.30 1.0) \lvec(0.70 1.0)
  \move(0.28 0.95) \lvec(0.72 0.95)
  \move(0.27 0.9) \lvec(0.73 0.9)
  \move(0.25 0.85) \lvec(0.75 0.85)
  \move(0.23 0.8) \lvec(0.77 0.8)
  \move(0.21 0.75) \lvec(0.79 0.75)
  \move(0.19 0.7) \lvec(0.81 0.7)
  \move(0.16 0.65) \lvec(0.40 0.65)
  \move(0.60 0.65) \lvec(0.84 0.65)
  \move(0.13 0.6) \lvec(0.30 0.6)
  \move(0.70 0.6) \lvec(0.87 0.6)
  \move(0.11 0.55) \lvec(0.25 0.55)
  \move(0.75 0.55) \lvec(0.89 0.55)
  \move(0.09 0.5) \lvec(0.20 0.5)
  \move(0.80 0.5) \lvec(0.91 0.5) \lpatt()
  \textref h:C v:C \htext(0.5 0.35){$\frac{1}{2}(1 + i)$}
  \move(0.5 0.7) \setgray 0.5 \linewd 0.04 \lvec(0.5 1.5)
  \move(1 0) \larc r:1 sd:108 ed:175
  \move(0 0) \larc r:1 sd:5 ed:72
  \setgray 0 \linewd 0.005
  \move(0.5 0.7) \fcir f:0 r:0.025
  \move(1 0.7) \arrowheadtype t:F \avec(0.55 0.7)
  \textref h:C v:C \htext(1.3 0.7){$\frac{1}{2} + b_1 i$}

\end{texdraw}
\end{center}
\caption{$\Omega_5$ contains a neighborhood of $e^{\pi i/3}$.}
\end{figure}

For general flat tori, a computer simulation suggests the
following picture: Let $\Omega_3$ (resp.\ $\Omega_5$) be the
subset of the moduli space $\mathcal{M}_1\cup \{\infty\}\cong S^2$
which corresponds to tori with three (resp.\ five) critical
points, then $\Omega_3 \cup\{\infty\}$ is a closed subset
containing $i$, $\Omega_5$ is an open subset containing $e^{\pi
i/3}$, both of them are simply connected and their common boundary
$C := \p\Omega_3 = \p\Omega_5$ is a curve homeomorphic to $S^1$
containing $\infty$. Moreover, the extra critical points are split
out from some half period point when the tori move from $\Omega_3$
to $\Omega_5$ across $C$.

Concerning with the experimental observation, we propose to prove
it by the method of deformations in $\mathcal{M}_1$. The
degeneracy analysis of critical points, especially the half period
points, is a crucial step. In this direction we have the following
partial result on tori corresponding to the line ${\rm Re}\, \tau
= 1/2$. These are equivalent to the rhombus tori and $\tau =
\frac{1}{2}(1 + i)$ is equivalent to the square torus where there
are only three critical points.

\begin{theorem}[Moduli Dependence]
Let $\omega_1 = 1$ and $\omega_2 = \tau = \frac{1}{2} + ib$ with
$b > 0$. Then
\begin{itemize}
\item [(1)] There exists $b_0 < \frac{1}{2} < b_1 < \sqrt 3/2$
such that $\frac{1}{2}{\omega_1}$ is a degenerate critical point
of $G(z;\tau)$ if and only if $b = b_0$ or $b = b_1$. Moreover,
$\frac{1}{2}{\omega_1}$ is a local minimum point of $G(z;\tau)$ if
$ b_0 < b < b_1$ and it is a saddle point if $b < b_0$ or $b
> b_1$.

\item [(2)] Both $\frac{1}{2}\omega_2$ and $\frac{1}{2}\omega_3$
are non-degenerate saddle points of $G$.

\item [(3)] $G(z;\tau)$ has two more critical points $\pm
z_0(\tau)$ when $b < b_0$ or $b > b_1$. They are non-degenerate
global minimum points of $G$ and in the former case
\begin{equation*}
{\rm Re}\,z_0(\tau) = \frac{1}{2};\quad  0< {\rm Im}\, z_0(\tau) <
\frac b 2.
\end{equation*}
\end{itemize}
\end{theorem}

Part (1) gives a strong support of the conjectural shape of the
decomposition $\mathcal{M}_1 = \Omega_3\cup \Omega_5$. Part (3)
implies that minimizers of $J_{8\pi}$ exist for tori with $\tau =
\frac{1}{2} + ib$ with $b < b_0$ or $b > b_1$.

The proofs are given in \S6. Notably Lemma 6.1, 6.2 and Theorem
6.6, 6.7. They rely on two fundamental inequalities on special
values of elliptic functions and we would like to single out the
statements (recall that $e_i = \wp(\frac{1}{2}\omega_i)$ and
$\eta_i = 2\zeta(\frac{1}{2}\omega_i)$):

\begin{theorem}[Fundamental Inequalities]\label{thm:1-3}
Let $\omega_1 = 1$ and $\omega_2 = \tau = \frac{1}{2} + ib$ with
$b > 0$. Then
\begin{itemize}
\item [(1)] $\displaystyle{\frac{d}{db}(e_1 + \eta_1) =
-4\pi\frac{d^2}{db^2}\log|\T_2(0)| > 0}$. \item [(2)]
$\displaystyle{\frac{1}{2} e_1 - \eta_1 =
4\pi\frac{d}{db}\log|\T_3(0)| < 0}$ and
$\displaystyle{\frac{d^2}{db^2}\log|\T_3(0)|
> 0}$. The same holds for $\T_4(0) = \overline{\T_3(0)}$.
In particular, $e_1$ increases in $b$.
\end{itemize}
\end{theorem}

\begin{figure}
\centering
\includegraphics[width=0.5\textwidth]{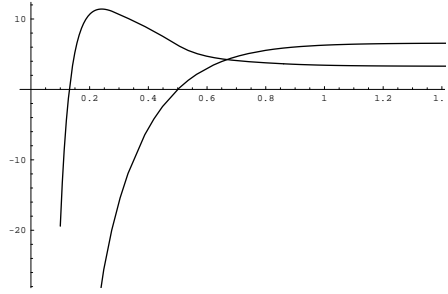}
\renewcommand{\figurename}{Figure}
\caption{Graphs of $\eta_1$ (the left one) and $e_1$ in $b$ where
$e_1$ is increasing.}
\end{figure}

\begin{figure}
\centering
\includegraphics[width=0.5\textwidth]{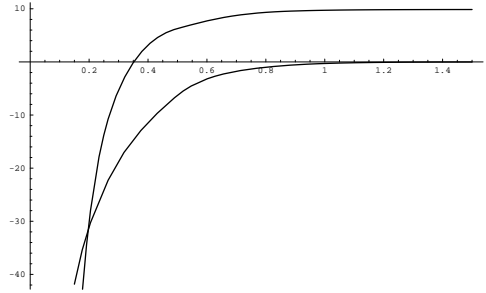}
\renewcommand{\figurename}{Figure}
\caption{Graphs of $e_1 + \eta_1$ (the upper one) and $\frac{1}{2}
e_1 - \eta_1$. Both functions are increasing in $b$ and
$\frac{1}{2} e_1 - \eta_1 \nearrow 0$.}
\end{figure}

These modular functions come into play due to the explicit
computation of Hessian at half periods along ${\rm Re}\,\tau =
\frac{1}{2}$ (cf.\ (\ref{eq:4-3}) and (\ref{eq:4-21})):
\begin{align*}
4\pi^2 \det D^2 G\Big(\frac{\omega_1}{2}\Big) &= (e_1 +
\eta_1)\Big(e_1 + \eta_1 - \frac{2\pi}{b}\Big),\\
4\pi^2 \det D^2 G\Big(\frac{\omega_2}{2}\Big) &= \Big(\eta_1 -
\frac 1 2 e_1\Big) \Big(\frac{2\pi}{b} + \frac 1 2 e_1 -
\eta_1\Big) - |{\rm Im}\, e_2|^2.
\end{align*}

Although $e_i$'s and $\eta_i$'s are classical objects, we were
unable to find an appropriate reference where these inequalities
were studied. Part of (2), namely $\frac{1}{2}e_1 - \eta_1 < 0$,
can be proved within the Weierstrass theory (cf.\
(\ref{eq:4-23})). The whole theorem, however, requires theta
functions in an essential way. Theta functions are recalled in \S7
and the theorem is proved in Theorem 8.1 and Theorem 9.1. The
proofs make use of the modularity of special values of theta
functions (Jacobi's imaginary transformation formula) as well as
the Jacobi triple product formula. Notice that the geometric
meaning of these two inequalities has not yet been fully explored.
For example, the variation on signs from $\vartheta_2$ to
$\vartheta_3$ is still mysterious to us.\medskip

{\bf Acknowledgements}: C.-S.~Lin is partially supported by the
National Science Council. C.-L.~Wang is partially supported by the
Shiing-shen Chern Fellowship of the National Center for Theoretic
Sciences at Hsinchu, Taiwan.

\section{Green Functions and Periods Integrals}
\setcounter{equation}{0}

We start with some basic properties of the Green functions that
will be used in the proof of Theorem 1.1. Detailed behavior of the
Green functions and their critical points will be studied in later
sections.

Let $T = \mathbb{C}/\Bbb Z\omega_1 + \Bbb Z\omega_2$ be a flat
torus. As usual we let $\omega_3 = \omega_1 + \omega_2$. The Green
function $G(z, w)$ is the unique function on $T$ which satisfies
\begin{equation}
-\triangle_z G(z, w) = \delta_{w}(z) - \frac{1}{|T|}
\end{equation}
and $\int_T G(z, w)\,dA = 0$. It has the property that $G(z, w) =
G(w, z)$ and it is smooth in $(z, w)$ except along the diagonal $z
= w$, where
\begin{equation}\label{eq:2-2}
G(z, w) = -\frac{1}{2\pi}\log|z - w| + O(|z - w|^2) + C
\end{equation}
for a constant $C$ which is independent of $z$ and $w$. Moreover,
due to the translation invariance of $T$ we have that $G(z, w) =
G(z - w, 0)$. Hence it is also customary to call $G(z) := G(z, 0)$
the Green function. It is an even function with the only
singularity at 0.

There are explicit formulae for $G(z, w)$ in terms of elliptic
functions, either in terms of the Weierstrass $\wp$ function or
the Jacobi-Riemann theta functions $\T_j$. Both are developed in
this paper since they have different advantages. We adopt the
first approach in this section.

\begin{lemma}\label{green-w}
There exists a constant $C(\tau)$, $\tau = \omega_2/\omega_1$,
such that
\begin{equation}{\label{eq:3-1}}
8\pi G(z) = \frac{2}{|T|}\int_T \log |\wp(\xi)-\wp(z)|\, dA +
C(\tau),
\end{equation}
\end{lemma}

It is straightforward to verify that the function of $z$, defined
in the right hand side of (\ref{eq:3-1}), satisfies the equation
for the Green function. By comparing with the behavior near 0 we
obtain Lemma \ref{green-w}. Since the proof is elementary, also an
equivalent form in theta functions will be proved in Lemma
\ref{green-theta}, we skip the details here.

In view of Lemma \ref{green-w}, in order to analyze critical
points of $G(z)$, it is natural to consider the following periods
integral
\begin{equation}{\label{eq:204}}
F(z)=\int_L \frac{\wp'(z)}{\wp(\xi) - \wp(z)}\, d\xi,
\end{equation}
where $L$ is a line segment in $T$ which is parallel to the
$\omega_1$-axis.

Fix a fundamental domain $T^0 = \{\,s\omega_1 + t\omega_2 \,|\,
-\frac{1}{2} \le s, t \le \frac{1}{2}\,\}$ and set $L^* = -L$.
Then $F(z)$ is an analytic function, except at 0, in each region
of $T^0$ divided by $L \cup L^*$. Clearly, ${\wp'(z)}/(\wp(\xi) -
\wp(z))$ has residue $\pm 1$ at $\xi = \pm z$. Thus for any fixed
$z$, $F(z)$ may change its value by $\pm 2\pi i$ if the
integration lines cross $z$. Let $T^0 \backslash (L \cup L^*) =
T_1\cup T_2\cup T_3$, where $T_1$ is the region above $L \cup
L^*$, $T_2$ is the region bounded by $L$ and $L^*$ and $T_3$ is
the region below $L \cup L^*$. Recall that $\zeta'(z)=-\wp(z)$ and
$\eta_i = \zeta(z + \omega_i) - \zeta(z)$ for $i\in\{1, 2\}$. Then
we have
\begin{lemma}\label{lem:2-2}
Let $C_1 = 2\pi i$, $C_2 = 0$ and $C_3 = -2\pi i$. Then for $z \in
T_k$,
\begin{equation}{\label{eq:205}}
F(z) = 2\omega_1\zeta(z) - 2\eta_1 z + C_k.
\end{equation}
\end{lemma}

\begin{proof} For $z \in T_1 \cup T_2 \cup T_3$, we have
$$
F'(z) = \int_L \frac{d}{dz} \left( \frac{\wp'(z)}{\wp(\xi) -
\wp(z)} \right) d\xi.
$$
Clearly, $z$ and $-z$ are the only (double) poles of
$\displaystyle{\frac{d}{dz} \left(\frac{\wp'(z)}{\wp(\xi) -
\wp(z)}\right)}$ as a meromorphic function of $\xi$ and
$\displaystyle{\frac{d}{dz} \left(\frac{\wp'(z)}{\wp(\xi) -
\wp(z)}\right)}$ has zero residues at $\xi = z$ and $-z$. Thus the
value of $F'(z)$ is independent of $L$ and it is easy to see
$F'(z)$ is a meromorphic function with the only singularity at 0.

By fixing $L$ such that $0\not\in L \cup L^*$, a straightforward
computation shows that
$$
F(z) = \frac{2\omega_1}{z} - 2\eta_1 z + O(z^2)
$$
in a neighborhood of $0$. Therefore
$$
F'(z) = -2\omega_1 \wp(z) - 2\eta_1.
$$

By integrating $F'$, we get
$$
F(z) = 2\omega_1 \zeta(z) - 2\eta_1 z + C_k.
$$

Since $F({\omega_2}/{2}) = 0, F({\omega_1}/{2}) = 0$ and
$F(-{\omega_2}/{2}) = 0$, $C_k$ is as claimed. Here we have used
the Legendre relation $\eta_1 \omega_2-\eta_2 \omega_1 = 2\pi i$.
\end{proof}

From (\ref{eq:3-1}), we have
\begin{equation}{\label{eq:3-2}}
8\pi G_z = \frac{1}{|T|}\int_T\;\frac{-\wp'(z)}{\wp(\xi) -
\wp(z)}\, dA.
\end{equation}

\begin{lemma}\label{critical}
Let $G$ be the Green function. Then for $z = t\omega_1 +
s\omega_2$,
\begin{equation}{\label{eq:3-3}}
G_z = -\frac{1}{4\pi}(\zeta(z) - \eta_1t - \eta_2 s).
\end{equation}
In particular, $z$ is a critical point of $G$ if and only if
\begin{equation}{\label{eq:3-4}}
\zeta(t\omega_1 + s\omega_2) = t\eta_1 + s\eta_2.
\end{equation}
\end{lemma}

\begin{proof} We shall prove (\ref{eq:3-3}) by applying Lemma
\ref{lem:2-2}. Since critical points appear in pair, without loss
of generality we may assume that $z = t\omega_1 + s\omega_2$ with
$s \ge 0$. We first integrate (\ref{eq:3-2}) along the $\omega_1$
direction and obtain
\begin{equation*}
\begin{split}
f(s_0) :&= \int_{L_1(s_0)} \frac{-\wp'(z)}{\wp(\xi)-\wp(z)}\,
d\xi\\
&= \left\{\begin{array}{ll}
-2\omega_1\zeta(z) + 2\eta_1 z & \mbox{if} \quad s_0 > s,\\
-2\omega_1\zeta(z) + 2\eta_1 z - 2\pi i\quad & \mbox{if}
\quad -s < s_0 < s,\\
-2\omega_1\zeta(z) + 2\eta_1 z & \mbox{if} \quad s_0 < -s,
\end{array}\right.
\end{split}
\end{equation*}
where $L_1(s_0) = \{\, t_0\omega_1 + s_0\omega_2 \mid |t_0|\leq
\frac 1 2\,\}$. Thus,
\begin{equation*}
\begin{split}
8\pi G_z &= \int_{-\frac 1 2}^{\frac 1 2} \int_{L_1(s_0)}
\frac{-\wp'(z)}{\wp(\xi) - \wp(z)}\,
dt_0\, ds_0 = \omega_1^{-1} \int_{-\frac 1 2}^{\frac 1 2} f(s_0)\, ds_0\\
&= \omega_1^{-1}((-2\omega_1\zeta(z) + 2\eta_1z) (1 - 2s) +
(-2\omega_1\zeta(z) + 2\eta_1 z - 2\pi i)2s)\\
&= \omega_1^{-1}(-2\omega_1\zeta(z) + 2\eta_1 z - 4\pi si)\\
&= \omega_1^{-1}(-2\omega_1\zeta(z) + 2\eta_1 t\omega_1 + 2s
(\eta_1 \omega_2 - 2\pi i))\\
&= -2\zeta(z) + 2\eta_1 t + 2s\eta_2,
\end{split}
\end{equation*}
where the Legendre relation is used again. \end{proof}

\begin{corollary}\label{lem:2-3}
Let $G(z)$ be the Green function. Then $\frac{1}{2}\omega_k$, $k
\in \{1, 2, 3\}$ are critical points of $G(z)$. Furthermore, if
$z$ is a critical point of $G$ then both periods integrals
\begin{equation}{\label{eq:206}}
\begin{split}
F_1 (z) &:= 2(\omega_1 \zeta (z) - \eta_1 z) \quad \mbox{and}\\
F_2 (z) &:= 2(\omega_2 \zeta (z) - \eta_2 z)
\end{split}
\end{equation}
are purely imaginary numbers.
\end{corollary}

\begin{proof} The half-periods $\frac{1}{2}{\omega_1}$,
$\frac{1}{2}{\omega_2}$ and $\frac{1}{2}{\omega_3}$ are obvious
solutions of (\ref{eq:3-4}). Alternatively, the half periods are
critical points of any even functions. Indeed for $G(z) = G(-z)$,
we get $\nabla G(z) = -\nabla G(-z)$. Let $p =
\frac{1}{2}{\omega_i}$ for some $i\in\{1,2,3\}$, then $p = -p$ in
$T$ and so $\nabla G(p) = -\nabla G(p) = 0$.

If $z = t\omega_1 + s\omega_2$ is a critical point, then by Lemma
\ref{critical},
$$
\omega_1\zeta(z) - \eta_1 z = \omega_1(t\eta_1 + s\eta_2) - \eta_1
(t\omega_1 + s\omega_2) = s(\omega_1\eta_2 - \omega_2 \eta_1) =
-2s\pi i.
$$

The proof for $F_2$ is similar.
\end{proof}

\begin{example}
When $\tau = \omega_2/\omega_1 \in i\,\mathbb{R}$, by symmetry
considerations it is known (cf.\ \cite{CLW}, Lemma 2.1) that the
half periods are all the critical points.
\end{example}

\begin{example}
There are tori such that equation (\ref{eq:3-4})
has more than three solutions. One such example is the torus with
$\omega_1 = 1$ and $\omega_2 = \frac 1 2 (1 + \sqrt{3}i)$. In this
case, the multiplication map $z\mapsto \omega_2 z$ is simply the
counterclockwise rotation by angle $\pi/3$, which preserves the
lattice $\mathbb{Z}\omega_1 + \mathbb{Z}\omega_2$, hence $\wp$
satisfies
\begin{equation}{\label{eq:305}}
\wp(\omega_2 z) = \wp(z)/\omega_2^2.
\end{equation}

Similarly in $\wp'^2 = 4\wp^3 - g_2 \wp - g_3$,
$$
g_2 = 60\sum\nolimits' \frac{1}{\omega^4} =
60\sum\nolimits'\frac{1}{(\omega_2\omega)^4} = \omega_2^2 g_2,
$$
which implies that $g_2 = 0$ and
\begin{equation}{\label{eq:3-5}}
\wp''=6\wp^2.
\end{equation}

Let $z_0$ be a zero of $\wp(z)$. Then $\wp''(z_0)=0$ too. By
(\ref{eq:305}), $\wp(\omega_2z_0)=0$, hence either $\omega_2
z_0=z_0$ or $\omega_2 z_0=-z_0$ on $T$ since $\wp(z)=0$ has zeros
at $z_0$ and $-z_0$ only. From here, it is easy to check that
either $z_0$ is one of the half periods or $z_0 = \pm
\frac{1}{3}\omega_3$. But $z_0$ can not be a half period because
$\wp''(z_0) \neq  0$ at any half period. Therefore, we conclude
that $z_0 = \pm\frac{1}{3}\omega_3$ and
$\wp''(\pm\frac{1}{3}\omega_3) = \wp(\pm\frac{1}{3}\omega_3) = 0$.

We claim that $\frac{1}{3}\omega_3$ is a critical point. Indeed
from the addition formula
\begin{equation}{\label{eq:3-6}}
\zeta(2z)=2\zeta(z)+\frac 1 2 \frac{\wp''(z)}{\wp'(z)}
\end{equation}
we have
\begin{equation}{\label{eq:3-7}}
\zeta\Big(\frac{2\omega_3}{3}\Big) =
2\zeta\Big(\frac{\omega_3}{3}\Big).
\end{equation}
On the other hand,
$$
\zeta\Big(\frac{2\omega_3}{3}\Big) =
\zeta\Big(\frac{-\omega_3}{3}\Big) + (\eta_1 + \eta_2) =
-\zeta\Big(\frac{\omega_3}{3}\Big) + \eta_1 + \eta_2.
$$ Together with
(\ref{eq:3-7}) we get
$$
\zeta\Big(\frac{\omega_3}{3}\Big) = \frac 1 3 (\eta_1+\eta_2).
$$
That is, $\frac{1}{3}\omega_3$ satisfies the critical point
equation.

Thus $G(z)$ has at least five critical points at
$\frac{1}{2}\omega_k$, $k=1,2,3$ and $\pm \frac{1}{3}\omega_3$
when $\tau = \omega_2/\omega_1 = \frac{1}{2}(1 + \sqrt{3}i)$.

By way of Theorem \ref{Green}, these are precisely the five
critical points, though we do not know how to prove this directly.
\end{example}

To conclude this section, let $u$ be a solution of (\ref{1-1})
with $\rho = 8\pi$ and set
\begin{equation}{\label{eq:2-6}}
v(z) = u(z) + 8\pi G(z).
\end{equation}
Then $v(z)$ satisfies
\begin{equation}{\label{eq:2-7}}
\triangle v(z) + 8\pi \Big(e^{-8\pi G(z)} e^{v(z)} - \frac{1}{|T|}
\Big) = 0
\end{equation}
in $T$. By (\ref{eq:2-2}), it is obvious that $v(z)$ is a smooth
solution of (\ref{eq:2-7}). An important fact which we need is the
following: Assume that there is a blow-up sequence of solutions
$v_j(z)$ of (\ref{eq:2-7}). That is,
$$
v_j(p_j) = \max_T v_j \to  +\infty \quad \mbox{as}\quad j\to
\infty.
$$
Then the limit $p = \lim_{j\to \infty} p_j$ is the only blow-up
point of $\{v_j\}$ and $p$ is in fact a critical point of $G(z)$:
\begin{equation}{\label{eq:2-8}}
\nabla G(p)=0.
\end{equation}
We refer the reader to \cite{Chen-Lin1} (p.\ 739, Estimate B) for
a proof of (\ref{eq:2-8}).

\section{The Criterion for Existence via Monodromies}
\setcounter{equation}{0}

Consider the mean field equation
\begin{equation}{\label{3-1}}
\triangle u + \rho e^u = \rho\delta_0,\quad \rho\in \mathbb{R}^+
\end{equation}
in a flat torus $T$, where $\delta_0$ is the Dirac measure with
singularity at 0 and the volume of $T$ is normalized to be 1. A
well known theorem due to Liouville says that any solution $u$ of
$\triangle u + \rho e^u = 0$ in a simply connected domain $\Omega
\subset \mathbb{C}$ must be of the form
\begin{equation}{\label{3-2}}
u = c_1 + \log \frac{|f'|^2}{(1 + |f|^2)^2},
\end{equation}
where $f$ is holomorphic in $\Omega$. Conventionally $f$ is called
a developing map of $u$. Given a torus $T = \mathbb{C}/\Bbb Z
\omega_1 + \Bbb Z\omega_2$, by gluing the $f$'s among simply
connected domains it was shown in \cite{CLW} that for $\rho = 4\pi
l$, $l \in \mathbb{N}$, (\ref{3-2}) holds on the whole
$\mathbb{C}$ with $f$ a meromorphic function. (The statement there
is for rectangular tori with $l = 2$, but the proof works for the
general case.)

It is straightforward to show that $u$ and $f$ satisfy
\begin{equation}{\label{3-3}}
u_{zz} - \frac{1}{2} u_z^2 = \frac{f'''}{f'} - \frac{3}{2}
\left(\frac{f''}{f'}\right)^2.
\end{equation}

The right hand side of (\ref{3-3}) is the Schwartz derivative of
$f$. Thus for any two developing maps $f$ and $\tilde f$ of $u$,
there exists $S = \begin{pmatrix} p & -\bar q\\ q & \bar
p\end{pmatrix}\in {\rm PSU}(1)$ (i.e. $p$, $q\in \Bbb C$ and
$|p|^2 + |q|^2 = 1$) such that
\begin{equation}{\label{3-4}}
\tilde f = Sf := \frac{pf - \bar q}{qf + \bar p}.
\end{equation}

Now we look for the constraints. The first type of constraints are
imposed by the double periodicity of the equation. By applying
(\ref{3-4}) to $f(z + \omega_1)$ and $f(z + \omega_2)$, we find
$S_1$ and $S_2$ in ${\rm PSU}(1)$ with
\begin{equation}{\label{can}}
\begin{split}
f(z + \omega_1) &= S_1 f,\\
f(z + \omega_2) &= S_2 f.
\end{split}
\end{equation}
These relations also force that $S_1S_2 = S_2 S_1$ (up to a sign,
as $A \equiv -A$ in ${\rm PSU}(1)$).

The second type of constraints are imposed by the Dirac
singularity of (\ref{3-1}) at 0. A straightforward local
computation with (\ref{3-2}) shows that

\begin{lemma}\label{order}
\begin{itemize}
\item[(1)] If $f(z)$ has a pole at $z_0 \equiv 0 \pmod{\omega_1,
\omega_2}$, then the order $k = l + 1$.

\item[(2)] If $f(z) = a_0 + a_k z^k + \cdots$ is regular at $z_0
\equiv 0 \pmod{\omega_1, \omega_2}$ with $a_k \ne 0$ then $k = l +
1$.

\item[(3)] If $f(z)$ has a pole at $z_0 \not\equiv 0
\pmod{\omega_1, \omega_2}$, then the order is 1.

\item[(4)] If $f(z) = a_0 + a_k (z - z_0)^k + \cdots$ is regular at
$z_0 \not\equiv 0 \pmod{\omega_1, \omega_2}$ with $a_k \ne 0$ then
$k = 1$.\smallskip
\end{itemize}
\end{lemma}

Now we are in a position to prove Theorem 1.1, namely the case $l
= 2$.

\begin{proof} We first prove the ``only if'' part. Let $u$ be a solution
and $f$ be a developing map of $u$. By the above discussion, we
may assume, after conjugating a matrix in ${\rm PSU}(1)$, that
$S_1 = \begin{pmatrix} e^{i\theta} & 0\\
0 & e^{-i\theta}\end{pmatrix}$ for some $\theta\in \Bbb R$. Let
$S_2 =
\begin{pmatrix} p & -\bar q\\ q & \bar p\end{pmatrix}$ and then
(\ref{can}) becomes
\begin{equation}
\begin{split}
f(z + \omega_1) &= e^{2i\theta} f(z),\\
f(z + \omega_2) &= S_2 f(z).
\end{split}
\end{equation}

Since $S_1 S_2 = S_2 S_1$, a direct computation shows that there
are three possibilities:
\begin{itemize}
\item[(1)] $p=0$ and $e^{i\theta}=\pm i$; \item[(2)] $q=0$; and
\item[(3)] $e^{i\theta}=\pm 1$.
\end{itemize}\medskip

{\bf Case (1).} By assumption we have
\begin{equation}{\label{case-1-1}}
\begin{split}
f(z + \omega_1) &= -f(z),\\
f(z + \omega_2) &= -\frac{(\bar q)^2}{f(z)}.
\end{split}
\end{equation}

For any $l \in \mathbb{N}$, the logarithmic derivative
\begin{equation*}
g = (\log f)' = \frac{f'}{f}
\end{equation*}
is a non-constant elliptic function on $\tilde T = \mathbb{C}/\Bbb
Z\omega_1 + \Bbb Z 2\omega_2$ which has a simple pole at each zero
or pole of $f$. By Lemma \ref{order}, if $f$ or $1/f$ is singular
at $z = 0$ then $g$ has no zero, which is not possible. So $f$
must be regular at $z = 0$ with $f(0) \ne 0$. Moreover $g$ has two
zeros of order $l$ at $0$ and $\omega_2$.

Let $\wp(z)$, $\zeta(z)$ and $\sigma(z) = \exp \int^z \zeta(w)\,dw
= z + \cdots$ be the Weierstrass elliptic functions on $\tilde T$.
Recall that $\sigma$ is odd with a simple zero at each lattice
point. Moreover, for $\tilde \omega_1 = \omega_1$, $\tilde
\omega_2 = 2\omega_2$ and $\tilde \omega_3 = \omega_1 +
2\omega_2$,
\begin{equation} \label{sigma-period}
\sigma(z \pm \tilde\omega_i) = -e^{\mp\eta_i(z \pm
\frac{1}{2}\tilde{\omega_i})}\sigma(z).
\end{equation}

Now $l = 2$. By the standard representation of elliptic functions,
\begin{equation}
g(z) = A \frac{\sigma^2(z) \sigma^2(z - \omega_2)}{\sigma(z - a)
\sigma(z - b) \sigma(z - c) \sigma(z - d)}
\end{equation}
with four distinct simple poles such that $a + b + c + d =
2\omega_2$. We will show that such a function $g(z)$ does not
exist.

By (\ref{case-1-1}), we have $g(z + \omega_2) = -g(z)$. Hence we
may assume that $c \equiv a + \omega_2$ and $d \equiv b +
\omega_2$ modulo $\omega_1, 2\omega_2$. Thus $a + b \equiv 0$
modulo $\frac{1}{2}\omega_1, \omega_2$ and we arrive at two
inequivalent cases. (i) $(a, b, c, d) = (a, -a, a + \omega_2, -a +
\omega_2)$. (ii) $(a, b, c, d) = (a, -a + \frac{1}{2}\omega_1, a +
\omega_2, -a - \frac{1}{2}\omega_1 + \omega_2)$. Using
(\ref{sigma-period}), it is easily checked that (i) leads to $g(z
+ \omega_2) = g(z)$ and (ii) leads to $g(z + \omega_2) = -g(z)$.
Hence we are left with case (ii).

The residues of $g$ at $a, b, c$ and $d$ are equal to $-Ar$,
$Ar'$, $Ar$ and $-Ar'$ respectively, where
$$
r = \frac{\sigma^2(a + \omega_2) \sigma^2(a)}{\sigma(\omega_2)
\sigma(2a - \frac{1}{2}\omega_1 + \omega_2) \sigma(2a +
\frac{1}{2}\omega_1)}
$$
and
$$
r' = \frac{\sigma^2(a - \frac{1}{2}\omega_1) \sigma^2(a -
\frac{1}{2}\omega_1 + \omega_2)}{\sigma(2a - \frac{1}{2}\omega_1)
\sigma(2a - \frac{1}{2}\omega_1 + \omega_2) \sigma(\omega_1 -
\omega_2)}.
$$

We claim that $Ar = \pm 1$ and $Ar' = \pm 1$: Since
\begin{equation*}
f(z) = \exp \int^z g(w)\,dw
\end{equation*}
is well-defined, by the residue theorem, we must have $Ar = m$ for
some $m \in \mathbb{Z}$. Moreover one of $a$, $b$ is a pole of
order $|m|$ of $f$ and then by Lemma \ref{order} we conclude that
$m = \pm 1$. Similarly $Ar' = \pm 1$.

In particular we must have $r/r' = \pm 1$. Using
(\ref{sigma-period}), this is equivalent to
\begin{equation} \label{r/r'}
\frac{\sigma^2(a + \omega_2) \sigma^2(a)}{\sigma^2(a +
\frac{1}{2}\omega_1) \sigma^2(a - \frac{1}{2}\omega_1 + \omega_2)}
= \mp e^{\eta_1(-\frac{1}{2}\omega_1 + \omega_2)}.
\end{equation}

To solve $a$ from this equation, we first recall that
$$
\wp(z) - \wp(y) = -\frac{\sigma(z + y) \sigma(z - y)}{\sigma^2(z)
\sigma^2(y)}.
$$
By substituting $y = \frac{1}{2}\tilde \omega_i$ into it and using
(\ref{sigma-period}), we get
\begin{equation} \label{p-e}
\wp(z) - e_i = \frac{\sigma^2(z + \frac{1}{2}\tilde
\omega_i)}{\sigma^2(z) \sigma^2(\frac{1}{2}\tilde\omega_i)}\,
e^{-\eta_i z}.
\end{equation}

With (\ref{p-e}), the ``$+$'' case in equation (\ref{r/r'})
simplifies to
$$
\wp\Big(a - \frac{\omega_1}{2} + \omega_2\Big) - e_1 = \wp(a) -
e_1.
$$
That is, $2a \equiv \frac{1}{2}\omega_1 - \omega_2$. But this
implies that $b \equiv c$, a contradiction.

Similarly, the ``$-$'' case in (\ref{r/r'}) simplifies to (using
the Legendre relation)
$$
\wp\Big(a + \frac{\omega_1}{2}\Big) - e_3 = \wp(a) - e_3.
$$
That is, $2a + \frac{1}{2}\omega_1 \equiv 0$. This leads to $c
\equiv d$, which is again a contradiction.

\medskip

{\bf Case (2).} In this case we have
\begin{equation}{\label{eq:2-14}}
\begin{split}
f(z + \omega_1) &= e^{2i\theta_1}f(z),
\\f(z + \omega_2) &= e^{2i\theta_2}f(z).
\end{split}
\end{equation}

The logarithmic derivative $g = (\log f)' = f'/f$ is now elliptic
on $T$ which has a simple pole at each zero or pole of $f$. As in
Case (1), it suffices to investigate the situation when $f$ is
regular at 0 and $f(0)\neq 0$. Since $g$ has 0 as its only zero
(of order 2), we have
\begin{equation}{\label{eq:2-23}}
g(z) = A \frac{\sigma^2(z)}{\sigma(z - z_0) \sigma (z + z_0)}
\end{equation}
where $\sigma(z)$ is the Weierstrass sigma function on $T$.
Without loss of generality we may assume that $f$ has a zero at
$z_0$ and a pole at $-z_0$. In particular $z_0 \ne -z_0$ in $T$
and we conclude that $z_0\neq {\omega_k}/{2}$ for any $k\in
\{1,2,3\}$.

Notice that if a meromorphic function $f$ satisfies
(\ref{eq:2-14}), then $e^\lambda f$ also satisfies (\ref{eq:2-14})
for any $\lambda\in \Bbb R$. Thus
\begin{equation}\label{3-11}
u_\lambda(z) = c_1 + \log
\frac{e^{2\lambda}|f'(z)|^2}{(1+e^{2\lambda}|f(z)|^2)^2}
\end{equation}
is a scaling family of solutions of (\ref{3-1}).

Clearly $u_\lambda(z)\to -\infty$ as $\lambda\to +\infty$ for any
$z$ such that $f(z)\neq 0$ and $u_\lambda(z_0)\to +\infty$ as
$\lambda\to +\infty$. Hence $z_0$ is the blow-up point and we have
by (\ref{eq:2-8}) that
$$\nabla G(z_0)=0.$$
Namely, it is a critical point other than the half periods.
\medskip

{\bf Case (3).} In this case we get that $S_1$ is the identity. So
by another conjugation in ${\rm PSU}(1)$ we may assume that $S_2$
is in diagonal form. But this case is then reduced to Case (2).
The proof of the ``only if'' part is completed.\medskip

Now we prove the ``if'' part. Suppose that there is a critical
point $z_0$ of $G(z)$ with $z_0\neq \frac{1}{2}{\omega_k}$ for any
$k \in\{1,2,3\}$.

For any closed curve $C$ such that $z_0$ and $-z_0\not\in C$, the
residue theorem implies that
\begin{equation}{\label{eq:2-29}}
\int_C \frac{\wp'(z_0)}{\wp(\xi) - \wp(z_0)}\, d\xi = 2\pi m i
\end{equation}
for some $m \in\{1, 0, -1\}$. Hence
\begin{equation}{\label{eq:2-30}}
f(z) := \exp \left(\int_0^z \frac{\wp'(z_0)}{\wp(\xi) -
\wp(z_0)}\,d\xi\right)
\end{equation}
is well-defined as a meromorphic function. Notice that $f$ is
non-constant since $\wp'(z_0) \ne 0$.

Let $L_1$ and $L_2$ be lines in $T$ which are parallel to the
$\omega_1$-axis and $\omega_2$-axis respectively and with $\pm z_0
\not\in L_1, L_2$. Then for $j = 1, 2$,
\begin{equation}{\label{eq:2-31}}
f(z + \omega_j) = f(z) \exp \left(\int_{L_j}
\frac{\wp'(z_0)}{\wp(\xi) - \wp(z_0)}\, d\xi \right).
\end{equation}

By Lemma \ref{lem:2-2},
$$
\int_{L_j} \frac{\wp'(z_0)}{\wp(\xi) - \wp(z_0)}\, d\xi = F_j
(z_0) + C_k
$$
for some $C_k \in \{2\pi i, 0, -2\pi i\}$. Also by Corollary
\ref{lem:2-3},
$$
F_j(z_0) = 2(\omega_j \zeta (z_0) - \eta_j z_0) = 2i\theta_j \in
i\,\mathbb{R}.
$$
Hence for $j = 1, 2$,
\begin{equation}{\label{eq:22-34}}
f(z + \omega_j) = f(z) e^{2i\theta_j}.
\end{equation}
holds. Set
$$
u_\lambda(x) = c_1 + \log \frac{e^{2\lambda}|f'(z)|^2} {(1 +
e^{2\lambda}|f(z)|^2)^2}.
$$
Then $u_\lambda(x)$ satisfies (\ref{3-1}) for any $\lambda\in\Bbb
R$ and $u_\lambda$ is doubly periodic by (\ref{eq:22-34}).
Therefore, solutions have been constructed and the proof of
Theorem \ref{1-1} is completed.
\end{proof}

A similar argument leads to
\begin{theorem}\label{4pi}
For $\rho = 4\pi$, there exists an unique solution of (\ref{3-1}).
\end{theorem}

\begin{proof} By the same procedure of the previous proof, there are
three cases to be discussed. For Case (2), the subcases that $f$
or $1/f$ is singular at $z = 0$ leads to contradiction as before.
For the subcase that $f$ is regular at $z = 0$ and $f(0) \ne 0$ we
see that $f(z)/f'(z)$ is an elliptic function with $0$ as its only
simple pole (since now $k - 1 = l = 1$). Hence Case (2) does not
occur. Similarly Case (3) is not possible.

Now we consider Case (1). By (\ref{case-1-1}), the function
\begin{equation*}
g = (\log f)' = \frac{f'}{f}
\end{equation*}
is elliptic on $T'' = \mathbb{C}/\Bbb Z\omega_1 + \Bbb Z
2\omega_2$. $g$ has a simple pole at each zero or pole of $f$. By
Lemma \ref{order}, if $f$ or $1/f$ is singular at $z = 0$ then $g$
has no zero and we get a contradiction. So $f$ is regular at $z =
0$, $f(0) \ne 0$ and $g$ has two simple zeros at $0$ and
$\omega_2$.

Let $\sigma(z) = \exp \int^z \zeta(w)\,dw = z + \cdots$ be the
Weierstrass sigma function on $T''$. $\sigma$ is odd with a simple
zero at each lattice point. Then
\begin{equation}
g(z) = A\frac{\sigma(z)\sigma(z - \omega_2)}{\sigma(z - a)\sigma(z
- b)}
\end{equation}
for some $a$, $b$ with $a + b = \omega_2$.

From (\ref{case-1-1}), we have $g(z + \omega_2) = -g(z)$. So $a +
\omega_2 = b \pmod {\omega_1, 2\omega_2}$. Since the
representation of $g$ in terms of sigma functions is unique up to
the lattice $\{\omega_1, 2\omega_2\}$, there is an unique solution
of $(a, b)$:
\begin{equation}
a = -\frac{\omega_1}{2};\quad b = \frac{\omega_1}{2} + \omega_2.
\end{equation}

Notice that the residues of $g$ at $a$ and $b$ are equal to $-Ar$
and $Ar$ respectively, where
\begin{equation*}
r = \frac{\sigma(\frac{1}{2}\omega_1) \sigma(\frac{1}{2}\omega_1 +
\omega_2)}{\sigma(\omega_1 + \omega_2)}.
\end{equation*}

We claim that $Ar = \pm 1$. Since
\begin{equation*}
f(z) = \exp \int^z g(w)\,dw
\end{equation*}
is well defined, by the residue theorem, we must have $Ar = m$ for
some $m \in \mathbb{Z}$. Moreover one of $a$, $b$ is a pole of
order $|m|$ and then by Lemma \ref{order} we conclude that $m =
\pm 1$.

Conversely, by picking up $a$, $b$ and $A = 1/r$ as above, $f(z)$
is a uniquely defined meromorphic function up to a factor $f(0)$.
There is an unique choice of $f(0)$ up to a norm one factor such
that $ c:= f(\omega_2)f(0)$ has $|c| = 1$. Thus by integrating
$g(z + \omega_2) = -g(z)$ we get $f(z + \omega_2) = c/f(z)$.

By integrating $g(z + \omega_1) = g(z)$ we get $f(z + \omega_1) =
c'f(z)$ where
$$
c' = \frac{f(\omega_1)}{f(0)} = \exp \int_0^{\omega_1} g(z)\,dz.
$$
To evaluate the period integral, notice that
$g(\frac{1}{2}\omega_1 + u) = -g(\frac{1}{2}\omega_1 - u)$. By
using the Cauchy principal value integral and the fact that the
residue of $g$ at $\frac{1}{2}\omega$ is $\pm 1$, we get
\begin{equation}
\int_0^{\omega_1} g(z)\,dz = \pm \frac{1}{2} \times 2\pi i = \pm
\pi i
\end{equation}
and so $c' = -1$.

Thus $f$ gives rise to a solution of (\ref{3-1}) for $\rho =
4\pi$. The developing map for the other choice $A = -1/r$ is
$1/f(z)$ which leads to the same solution. The proof is completed.
\end{proof}

Since equation (\ref{3-1}) is invariant under $z \mapsto -z$, the
unique solution is necessarily even.

\section{An Uniqueness Theorem for $\rho \in [4\pi, 8\pi]$ via
Symmetrizations}

\setcounter{equation}{0}

From the previous section, for $\rho = 8\pi$, solutions to the
mean field equation exist in a one parameter scaling family in
$\lambda$ with developing map $f$ and centered at a critical point
other than the half periods. By choosing $\lambda = -\log |f(0)|$
we may assume that $f(0) = 1$. Then we have
$$
f(-z) = \frac{1}{f(z)}.
$$
Consider the particular solution
$$
u(z) = c_1 + \log \frac{|f'(z)|^2} {(1 + |f(z)|^2)^2}.
$$
It is easy to verify that $u(-z) = u(z)$ and $u$ is the unique
even function in this family of solutions. In order to prove the
uniqueness up to scaling, it is equivalent to prove the uniqueness
within the class of even functions.

The idea is to consider the following equation
\begin{equation}{\label{eq:1}}
\left\{
\begin{array}{l}
\displaystyle \triangle u + \rho e^u = \rho\delta_0
\quad \mbox{and}\\
u(-z) = u(z) \quad \mbox{on $T$}
\end{array}
\right.
\end{equation}
where $\rho\in [4\pi, 8\pi]$. We will use the {\it method of
symmetrization} to prove

\begin{theorem}\label{nondeg}
For $\rho \in [4\pi, 8\pi]$, the linearized equation of
(\ref{eq:1}) is non-degenerate. That is, the linearized equation
has only trivial solutions.
\end{theorem}

Together with the uniqueness of solution in the case $\rho = 4\pi$
(Theorem \ref{4pi}), we conclude the proof of Theorem \ref{uniq}
by the inverse function theorem.\medskip

We first prove Theorem \ref{sym}, the Symmetrization Lemma. The
proof will consist of several Lemmas. The first step is an
extension of the classical isoperimetric inequality of Bol for
domains in $\Bbb R^2$ with metric $e^w |dx|^2$ to the case when
the metric becomes singular.

Let $\Omega \subset \mathbb{R}^2$ be a domain and $w \in
C^2(\Omega)$ satisfy
\begin{equation}{\label{eq:2}}
\left\{
\begin{array}{ll}
\triangle w + e^w \geq 0 \quad\mbox{in $\Omega$ and}\\
\displaystyle\int_\Omega e^w\,dx \leq 8\pi.
\end{array}
\right.
\end{equation}
This is equivalent to saying that the Gaussian curvature of
$e^w|dx|^2$ is
$$
K \equiv -\frac{1}{2} e^{-w}{\triangle w} \le \frac 1 2.
$$

For any domain $\omega \Subset \Omega$, we set
\begin{equation}{\label{eq:3}}
m(\omega) = \int_\omega e^w\,dx \quad\mbox{and} \quad \ell(\p
\omega) = \int_{\p\omega} e^{w/2}\,ds.
\end{equation}
Bol's isoperimetric inequality says that if $\Omega$ is
simply-connected then
\begin{equation}{\label{eq:4}}
2\ell^2(\p\omega) \geq m(\omega)(8\pi - m(\omega)).
\end{equation}

We first extend it to the case when $w$ acquires singularities:
\begin{eqnarray}{\label{eq:5}}
\left\{
\begin{array}{ll}
\triangle w + e^w = \sum_{j = 1}^N 2\pi\alpha_j \delta_{p_j} \quad
\mbox{in $\Omega$ and}\\
\displaystyle\int_\Omega e^w\, dx \leq 8\pi,
\end{array}
\right.
\end{eqnarray}
with $\alpha_j > 0$, $j = 1,2,\ldots,N$.

\begin{lemma}\label{lem 4-2}
Let $\Omega$ be a simply-connected domain and $w$ be a solution of
(\ref{eq:5}), or more generally a sub-solution with prescribed
singularities: $w(x) - \sum_j \alpha_j \log |x - p_j| \in
C^2(\Omega)$. Then for any domain $\omega \Subset \Omega$, we have
\begin{equation}{\label{eq:6}}
2\ell^2(\p\omega)\geq m(\omega)(8\pi - m(\omega)).
\end{equation}
\end{lemma}

\begin{proof} Define $v$ and $w_\epsilon$ by
$$
w(x) = \sum_j \alpha_j \log |x - p_j| + v(x),
$$
$$
w_\varepsilon(x) = \sum_j \frac{\alpha_j}{2}\log (|x - p_j|^2 +
\varepsilon^2) + v(x).
$$
By straightforward computations, we have
\begin{align*}
&\triangle w_\varepsilon(x) + e^{w_\varepsilon(x)}\\
&= \sum_j \frac{2\alpha_j\varepsilon^2}{(|x - p_j|^2 +
\varepsilon^2)^2} + e^v \left( \prod_j (|x - p_j|^2 +
\varepsilon^2)^{{\alpha_j}/{2}} - \prod_j |x - p_j|^{\alpha_j}
\right) \geq 0.
\end{align*}
Let $\ell_\varepsilon$ and $m_\varepsilon$ be defined as in
(\ref{eq:3}) with respect to the metric
$e^{w_\varepsilon(x)}|dx|^2$. Then we have
$$
2\ell_\varepsilon^2 (\p\omega) \geq m_\varepsilon(\omega)(8\pi -
m_\varepsilon(\omega)) .
$$
By letting $\varepsilon\to 0$ we obtain (\ref{eq:6}).
\end{proof}

Next we consider the case that some of the $\alpha_j$'s are
negative. For our purpose, it suffices to consider the case with
only one singularity $p_1$ with negative $\alpha_1$ (and we only
need the case that $\alpha_1 = -1$). In view of (the proof of)
Lemma \ref{lem 4-2}, the problem is reduced to the case with only
one singularity $p_1$. In other words, let $w$ satisfy
\begin{equation}{\label{eq:7}}
\triangle w + e^w = -2\pi \delta_{p_1} \quad\mbox{in $\Omega$}.
\end{equation}

\begin{lemma}\label{isoperimetric}
Let $w$ satisfy (\ref{eq:7}) with $\Omega$ simply-connected.
Suppose that
\begin{equation}{\label{eq:8}}
\int_\Omega e^w dx\leq 4\pi,
\end{equation}
then
\begin{equation} \label{isoperi}
2\ell^2 (\p\omega) \geq m(\omega) (4\pi-m(\omega)).
\end{equation}
\end{lemma}

\begin{proof}
 We may assume that $p_1 = 0$. If $0\not\in \omega$ then
$$
2\ell^2 (\p\omega) \geq m (\omega) (8\pi - m(\omega)) >
m(\omega)(4\pi - m(\omega))
$$
by Bol's inequality trivially. If $0 \in \omega$, we consider the
double cover $\tilde \Omega$ of $\Omega$ branched at $0$. Namely
we set $\tilde\Omega = f^{-1}(\Omega)$ where
$$
x = f(z) = z^2 \quad \mbox{for $z\in\Bbb C$}.
$$
The induced metric $e^v|dz|^2$ on $\tilde \Omega$ satisfies
$$
e^{v(z)}|dz|^2 = e^{w(x)}|dx|^2 = e^{w(z^2)}4|z|^2|dz|^2.
$$
That is, the metric potential $v$ is the regular part
$$
v(z) := w(x) + \log |x| + \log 4 = w(z^2) + 2\log |z| + \log 4.
$$
By construction, $v$ satisfies
$$
\triangle v + e^v = 0\quad \mbox{in $\tilde \Omega\backslash
\{0\}$}.
$$
Since $v$ is bounded in a neighborhood of 0, by the regularity of
elliptic equations, $v(z)$ is smooth at 0. Hence $v$ satisfies
\begin{equation*} 
\triangle v + e^v = 0\quad \mbox{in $\tilde \Omega$}.
\end{equation*}

Let $\tilde \omega = f^{-1}(\omega)$. Clearly
$\p\tilde\omega\subseteq f^{-1}(\p\omega)$. Also
$$
l(\p\tilde \omega)\leq 2l(\p \omega) \quad \mbox{and} \quad
m(\tilde w)=2m(\omega),
$$
where
$$
l(\p\tilde\omega) =\int_{\p\tilde\omega} e^{v/2}\, ds
\quad\mbox{and}\quad m(\tilde \omega) =\int_{\tilde \omega} e^v\,
dx.
$$
By Bol's inequality, we have
$$
4\ell^2(\p\omega) \geq \ell^2(\p\tilde\omega) \geq \frac 1 2
m(\tilde \omega) (8\pi - m(\tilde \omega)) = m(\omega) (8\pi -
2m(\omega)).
$$
Thus $2 \ell^2 (\p\omega) \geq m(\omega) (4\pi - m(\omega))$.
\end{proof}

\begin{lemma}[Symmetrization I]\label{symmetry}
Let $\Omega \subset \mathbb{R}^2$ be a simply-connected domain
with $0\in\Omega$ and let $v$ be a solution of
$$
\triangle v + e^v = -2 \pi \delta_0
$$
in $\Omega$. If the first eigenvalue of $\triangle + e^v$ is zero
on $\Omega$ then
\begin{equation}{\label{eq:15}}
\int_\Omega  e^v\, dx \geq 2\pi.
\end{equation}
\end{lemma}

\begin{proof} Let $\psi$ be the first eigenfunction of $\triangle + e^v$:
\begin{equation}{\label{eq:16}}
\left\{
\begin{array}{l}
\triangle \psi + e^v \psi=0 \quad \mbox{in $\Omega$ and}\\
\psi = 0 \quad \mbox{on $\p\Omega$}.
\end{array}\right.
\end{equation}
In $\mathbb{R}^2$, let $U$ and $\varphi$ be the radially symmetric
functions
\begin{equation}{\label{eq:17}}
\left\{
\begin{array}{l}
\displaystyle U(x) = \log \frac{2}{(1+|x|)^2|x|} \quad \mbox{and}\\
\displaystyle \varphi(x)=\frac{1-|x|}{1+|x|}. \end{array}\right.
\end{equation}
From $\displaystyle{\triangle = \frac {\p^2}{\p r^2} + \frac {1}
{r} \frac{\p}{\p r} + \frac 1 {r^2}\frac{\p^2}{\p\theta^2}}$, it
is easy to verify that $U$ and $\varphi$ satisfy
\begin{equation}{\label{eq:18}}
\left\{
\begin{array}{l}
\triangle U + e^U = -2\pi\delta_0\quad \mbox{and}\\
\triangle\varphi +  e^U \varphi = 0\quad \mbox{in $\mathbb{R}^2$}.
\end{array}\right.
\end{equation}
For any $t > 0$, set $\Omega_t = \{ x\in\Omega \mid \psi(x)
> t\}$ and $r(t) > 0$ such that
\begin{equation}{\label{eq:19}}
\int_{B_{r(t)}} e^{U(x)}\, dx = \int_{\{\psi>t\}} e^{v(x)}\,dx,
\end{equation}
where $B_{r(t)}$ is the ball with center 0 and radius $r(t)$.
Clearly $r(t)$ is strictly decreasing in $t$ for $t\in (0, \max
\psi)$. In fact, $r(t)$ is Lipschitz in $t$. Denote by $\psi^*(r)$
the symmetrization of $\psi$ with respect to the measure
$e^{U(x)}\,dx$ and $e^{v(x)}\,dx$. That is,
$$
\psi^*(r) = \sup\{t\mid r < r(t)\}.
$$
Obviously $\psi^*(r)$ is decreasing in $r$ and for $t\in
(0,\max\psi)$, $\psi^*(r) = t$ if and only if $r(t) = r$. Thus by
(\ref{eq:19}), we have a decreasing function
\begin{equation}{\label{eq:20}}
f(t) := \int_{\{\psi^*>t\}}\, e^{U(x)}\, dx = \int_{\{\psi>t\}}
e^{v(x)}\, dx.
\end{equation}

By Lemma \ref{isoperimetric}, for any $t > 0$,
\begin{equation}{\label{eq:21}}
2\ell^2(\{\psi=t\}) \geq f(t)(4\pi - f(t)).
\end{equation}
We will use inequality (\ref{eq:21}) in the following computation:
For any $t>0$, by the Co-Area formula,
\begin{equation*}
\begin{split}
-\frac{d}{dt} \int_{\Omega_t} |\nabla\psi|^2\,dx =
\int_{\p\Omega_t} |\nabla \psi|\,ds \quad \mbox{and}\\
-f'(t) = -\frac{d}{dt} \int_{\Omega_t} e^{v}\,dx =
\int_{\p\Omega_t} \frac{e^v}{|\nabla \psi|}\,ds
\end{split}
\end{equation*}
hold almost everywhere in $t$. Thus
\begin{equation}{\label{eq:22}}
\begin{split}
& -\frac{d}{dt} \int_{\Omega_t} |\nabla\psi|^2\, dx =
\int_{\{\psi = t\}} |\nabla\psi|\, ds\\
& \geq \Biggl (\int_{\{\psi = t\}}  e^{v/2}\, ds\Biggl)^2 \Biggl(
\int_{\{\psi = t\}} \frac{e^{v}}{|\nabla \psi|}
\,ds\Biggl)^{-1}\\
&= -\ell^2(\{\psi = t\}) f'(t)^{-1}\\
& \geq -\frac 1 2 f(t)(4\pi - f(t)) f'(t)^{-1}.
\end{split}
\end{equation}
It is known that $\psi^* \in H^1_0(B_{r(0)})$ and the same
procedure for $\psi^*$ leads to
\begin{equation}
\begin{split}
&-\frac{d}{dt} \int_{\{\psi^* > t\}} |\nabla \psi^*|^2\, dx
= \int_{\{\psi^* = t\}} |\nabla \psi^*|\,ds\\
&= \Biggl(\int_{\{\psi^* = t\}} e^{U/2}\, ds\Biggl)^2 \Biggl(
\int_{\{\psi^* = t\}} \frac{e^{U}}{|\nabla \psi^*|}\,
ds\Biggl)^{-1}\\
&= -\ell^2(\{\psi^* = t\}) f'(t)^{-1} \\
&= -\frac 1 2 f(t)(4\pi - f(t)) f'(t)^{-1},
\end{split}
\end{equation}
with all inequalities being equalities. Hence
$$
-\frac{d}{dt} \int_{\{\psi>t\}} |\nabla \psi|^2\, dx \geq
-\frac{d}{dt} \int_{\{\psi^*>t\}} |\nabla\psi^*|^2ds
$$
holds almost everywhere in $t$. By integrating along $t$, we get
$$
\int_\Omega |\nabla\psi|^2 \,dx \geq \int_{B_{r(0)}}
|\nabla\psi^*|^2 \,dx.
$$
Since $\psi$ and $\psi^*$ have the same distribution (or by
looking at $-\int f'(t)t^2\,dt$ directly), we have
$$
\int_{B_{r(0)}} e^{U(x)} {\psi^*}^2\, dx = \int_\Omega
e^{v(x)}\psi^2\, dx.
$$
Therefore
$$
0 = \int_\Omega |\nabla\psi|^2\,dx - \int_\Omega e^v \psi^2\, dx
\geq \int_{B_{r(0)}} |\nabla\psi^*|^2\, dx - \int_{B_{r(0)}}
e^{U}{\psi^*}^2 \,dx.
$$
This implies that the linearized equation $\triangle + e^{U(x)}$
has non-positive first eigenvalue. By (\ref{eq:17}), this happens
if and only if $r(0) \geq 1$. Thus
$$
\int_\Omega e^{v(x)}\, dx =\int_{B_{r(0)}} e^{U(x)}\, dx \geq
2\pi.
$$
\end{proof}

\begin{remark}\label{easy-sym}
(See \cite{CCL}.) By applying symmetrization to $\triangle v + e^v
= 0$ in $\Omega$ with $\lambda_1(\triangle + e^v) = 0$, the
corresponding radially symmetric functions are
$$
U(x) = -2 \log\left(1 + \frac{1}{8}|x|^2\right) \quad
\mbox{and}\quad \varphi(x) = \frac{8 - |x|^2}{8 + |x|^2}.
$$
The same computations leads to $\int_\Omega e^v\,dx \ge 4\pi$.
\end{remark}

A closer look at the proof of Lemma \ref{symmetry} shows that it
works for more general situations as long as the isoperimetric
inequality holds:

\begin{lemma}[Symmetrization II]\label{symmetry+}
Let $\Omega \subset \mathbb{R}^2$ be a simply-connected domain and
let $v$ be a solution of
$$
\triangle v + e^v = \sum\nolimits_{j = 1}^N
2\pi\alpha_j\delta_{p_j}
$$
in $\Omega$. Suppose that the first eigenvalue of $\triangle +
e^v$ is zero on $\Omega$ with $\varphi$ the first eigenfunction.
If the isoperimetric inequality with respect to $ds^2 =
e^v|dx|^2$:
\begin{equation*}
2 \ell^2 (\p\omega) \geq m(\omega) (4\pi - m(\omega))
\end{equation*}
holds for all level domains $\omega = \{\varphi \ge t\}$ with $t
\ge 0$, then
\begin{equation*}{\label{eq:15+}}
\int_\Omega  e^v\, dx \geq 2\pi.
\end{equation*}
\end{lemma}

Notice that we do not need any further constraint on the sign of
$\alpha_j$.

For the last statement of Theorem \ref{sym}, the limiting
procedure of Lemma \ref{lem 4-2} implies that the isoperimetric
inequality (Lemma \ref{isoperimetric}) and symmetrization I (Lemma
\ref{symmetry}) both hold regardless on the presence of
singularities with non-negative $\alpha$.

Indeed the proof of Lemma \ref{isoperimetric} can be adapted to
the case
$$
\triangle w + e^w = -2\pi \delta_{p_1} + \sum\nolimits_{j = 2}^N
2\pi\alpha_j \delta_{p_j}
$$
with $\alpha_j > 0$ for $j = 2, \ldots, N$. On the double cover
$\bar \Omega \to \Omega$ branched over $p_1 = 0$, the metric
potential $v(z)$ again extends smoothly over $z = 0$ and satisfies
$$
\triangle v + e^v = \sum\nolimits_{j = 2}^N 2\pi\alpha_j
(\delta_{q_j} + \delta_{q'_j}),
$$
where $q_j, q'_j \in \bar \Omega$ are points lying over $p_j$. The
remaining argument works by using Lemma \ref{lem 4-2} and we still
conclude $2 \ell^2 (\p\omega) \geq m(\omega) (4\pi - m(\omega))$.

Thus the proof of Theorem \ref{sym} is completed.

\begin{proof} (of Theorem \ref{nondeg}.) Let $u$ be a solution of
equation (\ref{eq:1}). It is clear that we must have $\int_T e^u =
1$. Suppose that $\varphi(x)$ is a non-trivial solution of the
linearized equation at $u$:
\begin{equation}{\label{eq:9}}
\left\{
\begin{array}{l}
\triangle\varphi + \rho e^u \varphi = 0 \quad\mbox{and}\\
\varphi(z) = \varphi(-z) \quad \mbox{in $T$}.
\end{array}
\right.
\end{equation}
We will derive from this a contradiction.

Since both $u$ and $\varphi$ are even functions, by using $x =
\wp(z)$ as two-fold covering map of $T$ onto $S^2 = \mathbb{C}
\cup \{\infty\}$, we may require that $\wp$ being an isometry:
$$
e^{u(z)}|dz|^2 = e^{v(x)}|dx|^2 = e^{v(x)}|\wp'(z)|^2|dz|^2.
$$
Namely we set
\begin{equation}{\label{eq:10}}
v(x) := u(z) - \log|\wp'(z)|^2 \quad\mbox{and}\quad \psi(x) :=
\varphi(z).
\end{equation}
There are four branch points on $\mathbb{C} \cup\{\infty\}$,
namely $p_0 = \wp(0) = \infty$ and $p_j = e_j := \wp(\omega_j/2)$
for $j = 1, 2, 3$. Since $\wp'(z)^2 = 4\prod_{j = 1}^3 (x - e_j)$,
by construction $v(x)$ and $\psi(x)$ then satisfy
\begin{equation}{\label{eq:11}}
\left\{
\begin{array}{l}
\displaystyle \triangle v + \rho e^v = \sum\nolimits_{j = 1}^3
(-2\pi)\delta_{p_j} \quad\mbox{and}\\
\triangle \psi + \rho e^v \psi = 0\quad\mbox{in $\mathbb{R}^2$}.
\end{array}
\right.
\end{equation}

To take care of the point at infinity, we use coordinate $y = 1/x$
or equivalently we consider $T \to S^2$ via $y = 1/\wp(z) \sim
z^2$. The isometry condition reads as
$$
e^{u(z)}|dz|^2 = e^{w(y)}|dy|^2 =
e^{w(y)}\frac{|\wp'(z)|^2}{|\wp(z)|^4}|dz|^2.
$$
Near $y = 0$ we get
$$
w(y) = u(z) - \log\frac{|\wp'(z)|^2}{|\wp(z)|^4} \sim
\left(\frac{\rho}{4\pi} - 1\right)\log|y|.
$$
Thus $\rho \ge 4\pi$ implies that $p_0$ is a singularity with
non-negative $\alpha_0$:
$$
\triangle w + \rho e^w = \alpha_0 \delta_0 + \sum\nolimits_{j =
1}^3 (-2 \pi)\delta_{1/p_j}.
$$

In dealing with equation (\ref{eq:1}) and the above resulting
equations, by replacing $u$ by $u + \log \rho$ etc., we may (and
will) replace the $\rho$ in the left hand side by 1 for
simplicity. The total measure on $T$ and $\Bbb R^2$ are then given
by
$$
\int_T e^u \, dz = \rho \le 8\pi \quad \mbox{and} \quad \int_{\Bbb
R^2} e^v \, dx = \frac \rho 2 \le 4\pi.
$$

The nodal line of $\psi$ decomposes $S^2$ into at least two
connected components and at least two of them are simply
connected. If there is a simply connected component $\Omega$ which
contains no $p_j$'s, then the symmetrization (Remark
\ref{easy-sym}) leads to
$$
\int_{\Omega} e^v\,dx \geq 4\pi,
$$
which is a contradiction because $\mathbb{R}^2 \backslash \Omega
\ne \emptyset$. If every simply connected component $\Omega_i$, $i
= 1,\ldots,m$, contains only one $p_j$, then Lemma \ref{symmetry}
implies that
$$
\int_{\Omega_i} e^v\,dx \geq 2\pi \quad \mbox{for $i = 1, \ldots,
m$}.
$$
The sum is at least $2m\pi$, which is again impossible unless $m =
2$ and $\mathbb{R}^2 = \overline{\Omega}_1 \cup
\overline{\Omega}_2$. So without lost of generality we are left
with one of the following two situations:
$$
\mathbb{R}^2\cup \{\infty\}\backslash \{\psi=0\} = \Omega_+ \cup
\Omega_-
$$
where
$$
\Omega_+ \subset \{x\mid \psi(x) > 0\} \quad \mbox{and}\quad
\Omega_- \supset \{x\mid \psi(x) < 0\}.
$$
Both $\Omega_+$ and $\Omega_-$ are simply-connected.
\begin{itemize}
\item[(1)] Either $\Omega_-$ contains $p_1, p_2$ and $p_3 \in
\Omega_+$ or

\item[(2)] $p_1 \in \Omega_-$, $p_3\in \Omega_+$ and $p_2 \in C =
\p\Omega_+ = \p\Omega_-$.
\end{itemize}

Assume that we are in case (1). By Lemma \ref{isoperimetric}, we
have on $\Omega_+$
\begin{equation}{\label{eq:12}}
2 \ell^2 (\{\psi = t\})\geq m(\{\psi\geq t\})(4\pi - m (\{\psi\geq
t\}))
\end{equation}
for $t\geq 0$.

We will show that the similar inequality
\begin{equation}{\label{keypoint}}
2 \ell^2 (\{\psi = t\})\geq m(\{\psi\leq t\})(4\pi - m (\{\psi\leq
t\}))
\end{equation}
holds on $\Omega_-$ for all $t \le 0$.

Let $t \le 0$ and $\omega$ be a component of $\{\psi \le t\}$.

If $\omega$ contains at most one point of $p_1$ and $p_2$, then
Lemma \ref{isoperimetric} implies that
\begin{equation*}
2\ell^2(\p\omega) \geq (4\pi -m(\omega)) m(\omega).
\end{equation*}

If $\omega$ contains both $p_1$ and $p_2$, then $\mathbb{R}^2 \cup
\{\infty\} \backslash \omega$ is simply connected which contains
$p_3$ only. Thus by Lemma \ref{isoperimetric}
\begin{equation}\label{key}
\begin{split}
2\ell^2 (\p\omega) & \geq (4\pi - m(\mathbb{R}^2
\backslash\omega))
m(\mathbb{R}^2 \backslash\omega)\\
& = (4\pi - \rho/2 + m(\omega)) (\rho/2 - m(\omega))\\
& = (4\pi - m(\omega))m(\omega) + (4\pi - \rho/2) (\rho/2 -
2m(\omega)).
\end{split}
\end{equation}
Since $\rho\leq 8\pi$ and $\int_{\Omega_+} e^v dx \geq 2\pi$, we
get
$$
\frac \rho 2 = \int_{\mathbb{R}^2} e^v \geq \int_{\Omega_+} e^v +
\int_\omega e^v \geq 2\pi + m(\omega) \geq 2m(\omega).
$$
Then again
\begin{equation*}
2\ell^2(\p\omega) \geq (4\pi -m(\omega)) m(\omega)
\end{equation*}
with equality holds only when $\rho = 8\pi$ and $m(\omega) =
2\pi$.

Now it is a simple observation that domains which satisfy the
isoperimetric inequality (\ref{isoperi}) have the addition
property. Indeed, if $2a^2 \ge (4\pi - m)m$ and $2b^2 \ge (4\pi -
n)n$, then
\begin{align*}
2(a + b)^2 &= 2a^2 + 4ab + 2b^2 > (4\pi - m)m + (4\pi - n)n\\
&= 4\pi(m + n) - (m + n)^2 + 2mn > (4\pi - (m + n))(m + n).
\end{align*}
Hence (\ref{keypoint}) holds for all $t \le 0$.

Now we can apply Lemma \ref{symmetry+} to $\Omega_-$ to conclude
that
$$
\int_{\Omega_-} e^v\, dx = 2\pi
$$
(which already leads to a contradiction if $\rho < 8\pi$) and the
equality
$$
2\ell^2(\{\psi = t\}) = (4\pi - m(\{\psi\le t\}))m(\{\psi\le t\})
$$
in (\ref{keypoint}) holds for all $t\in (\min \psi, 0)$. This
implies that $\{\psi\le t\}$ has only one component and it
contains $p_1$ and $p_2$ for all $t\in (\min \psi, 0)$. But then
$\psi$ attains its minimum along a connected set $\psi^{-1}(\min
\psi)$ containing $p_1$ and $p_2$, which is impossible.

In case (2), $\rho < 8\pi$ again leads to a contradiction via the
same argument. For $\rho = 8\pi$, we have
$$
\int_{\Omega_+} e^v\, dx = \int_{\Omega_-} e^v\, dx = 2\pi
$$
and all inequalities in (\ref{eq:22}) are equalities. So under the
notations there
$$
\Biggl( \int_{\p\Omega_t}  e^{v/2}\, dx\Biggl)^2 =
\int_{\p\Omega_t} |\nabla \psi|\,ds \int_{\p\Omega_t}
\frac{e^v}{|\nabla\psi|}\,ds
$$
for all $t \in (\min_{\Omega_-} \psi, \max_{\Omega_+} \psi)$. This
implies that
$$
|\nabla\psi|^2 (x) = C(t) e^{v(x)}
$$
almost everywhere in $\psi^{-1}(t)$ for some constant $C$ which
depends only on $t$. By continuity we have
\begin{equation}\label{final}
|\nabla\psi|^2(x) = C(\psi(x)) e^{v(x)}
\end{equation}
for all $x$ except when $\psi(x) = \max_{\Omega_+} \psi$ or
$\psi(x)=\min_{\Omega_-}\psi$.

By letting $x = p_2\in \psi^{-1}(0)$, we find
$$
C(0) = |\nabla\psi|^2 (p_2) e^{-v(p_2)} = 0.
$$
By (\ref{final}) this implies that $|\nabla\psi(x)| = 0$ for all
$x \in \psi^{-1}(0)$, which is clearly impossible. Hence the proof
of Theorem \ref{nondeg} is completed.
\end{proof}

Since equation (\ref{eq:1}) has an unique solution at $\rho =
4\pi$, by the continuation from $\rho = 4\pi$ to $8\pi$ and
Theorem \ref{nondeg}, we conclude that (\ref{eq:1}) has a most a
solution at $\rho = 8\pi$, and it implies that the mean field
equation (\ref{1-1}) has at most one solution up to scaling. Thus,
Theorem \ref{uniq} is proved and then Theorem \ref{Green} follows
immediately.

\section{Comparing Critical Values of Green Functions}
\setcounter{equation}{0}

For simplicity, from now on we normalize all tori to have
$\omega_1 = 1$, $\omega_2 = \tau$.

In section 3 we have shown that the existence of solutions of
equation (\ref{1-1}) is equivalent to the existence of
non-half-period critical points of $G(z)$. The main goal of this
and next sections is to provide criteria for detecting minimum
points of $G(z)$. The following theorem is useful in this regard.

\begin{theorem}\label{thm:4-2}
Let $z_0$ and $z_1$ be two half-periods. Then $G(z_0)\geq G(z_1)$
if and only if $|\wp(z_0)|\geq |\wp(z_1)|$.
\end{theorem}

\begin{proof} By integrating (\ref{eq:3-3}), the Green function $G(z)$
can be represented by
\begin{equation}{\label{eq:3-9}}
G(z) = -\frac{1}{2\pi}{\rm Re}\int (\zeta(z) - \eta_1 z)\, dz +
\frac{1}{2b}y^2 + C(\tau).
\end{equation}
Thus
\begin{equation}{\label{eq:3-10}}
G\Big(\frac{\omega_2}{2}\Big) - G\Big(\frac{\omega_3}{2}\Big) =
\frac{1}{2\pi}{\rm Re}
\int_{\frac{\omega_2}{2}}^{\frac{\omega_3}{2}} (\zeta(z) -
\eta_1z)\,dz.
\end{equation}

Set $F(z) = \zeta(z) - \eta_1 z$. We have $F(z + \omega_1) = F(z)$
and
\begin{equation*}
\begin{split}
\zeta\Big(t + \frac{\omega_2}{2}\Big) - \eta_1 \Big(t +
\frac{\omega_2}{2}\Big) &= -\zeta\Big(\frac{-\omega_2}{2} - t\Big)
- \eta_1\Big(t +
\frac{\omega_2}{2}\Big)\\
&= -\zeta\Big(\frac{\omega_2}{2} - t\Big) + \eta_2 - \eta_1
\Big(t + \frac{\omega_2}{2}\Big)\\
&= -\Big[\zeta \Big(\frac{\omega_2}{2} - t\Big)-\eta_1
\Big(\frac{\omega_2}{2} - t\Big)\Big] + \eta_2
- \eta_1 \omega_2\\
&= -\Big[\zeta\Big(\frac{\omega_2}{2} - t\Big) - \eta_1
\Big(\frac{\omega_2}{2} - t\Big)\Big] - 2\pi i,
\end{split}
\end{equation*}
hence ${\rm Re}\,F(\frac{1}{2}\omega_2 + t)$ is antisymmetric in
$t\in \Bbb C$.

To calculate the integral in (\ref{eq:3-10}), we use the addition
theorem to get
\begin{equation*}
\begin{split}
\frac{\wp'(z)}{\wp(z) - e_1} &= \zeta\Big(z -
\frac{\omega_1}{2}\Big) + \zeta\Big(z +
\frac{\omega_1}{2}\Big) - 2\zeta(z)\\
&= F\Big(z - \frac{\omega_1}{2}\Big) + F\Big(z +
\frac{\omega_1}{2}\Big) - 2F(z).
\end{split}
\end{equation*}
Integrating it along the segment from $\frac{1}{2}\omega_2$ to
$\frac{1}{2}\omega_3$, we get
$$
\log \frac{e_3 - e_1}{e_2 - e_1} = \int_{L_1} F(z)\,dz +
\int_{L_2} F(z)\, dz - 2 \int_{L_3} F(z)\,dz,
$$
where $L_1$ is the line from $\frac{1}{2}(\omega_2 - \omega_1)$ to
$\frac{1}{2}\omega_2$, $L_2$ is the line from
$\frac{1}{2}\omega_3$ to $\frac{1}{2}\omega_2 + \omega_1$ and
$L_3$ is the line from $\frac{1}{2}\omega_2$ to
$\frac{1}{2}\omega_3$. Since $F(z) = F(z + \omega_1)$ and ${\rm
Re}\, F(z)$ is antisymmetric with respect to
$\frac{1}{2}\omega_2$, we have
$$
\log \frac{e_3 - e_1}{e_2 - e_1} = 2\int_L F(z)\, dz - 4
\int_{L_3} F(z)\, dz = -2\pi i - 4\int_{L_3} F(z)\, dz,
$$
where $L$ is the line from $\frac{1}{2}(\omega_2 - \omega_1)$ to
$\frac{1}{2}\omega_3$. Thus we have
$$
\log \left|\frac{e_3-e_1}{e_2-e_1}\right| = -4\,{\rm Re}\,
\int_{\frac{\omega_2}{2}}^{\frac{\omega_3}{2}} F(z)\,dz = -8\pi
\left(G\Big(\frac{\omega_2}{2}\Big) -
G\Big(\frac{\omega_3}{2}\Big)\right).
$$
That is,
\begin{equation}{\label{eq:3-11}}
G\Big(\frac{\omega_2}{2}\Big) - G\Big(\frac{\omega_3}{2}\Big) =
\frac{1}{8\pi} \log \left| \frac{e_2 - e_1}{e_3 - e_1}\right|.
\end{equation}

Similarly, by integrating (\ref{eq:3-3}) in the $\omega_2$
direction, we get
\begin{equation}
G\Big(\frac{\omega_1}{2}\Big) - G\Big(\frac{\omega_3}{2}\Big) =
\frac{1}{2\pi}{\rm Re}
\int_{\frac{\omega_1}{2}}^{\frac{\omega_3}{2}} (\zeta(z) - \eta_2
z)\,dz.
\end{equation}
The same proof then gives rise to
\begin{equation}\label{eq:5-5}
G\Big(\frac{\omega_1}{2}\Big) - G\Big(\frac{\omega_3}{2}\Big) =
\frac{1}{8\pi} \log \left|\frac{e_1 - e_2}{e_3 - e_2}\right|.
\end{equation}

By combining the above two formulae we get also that
\begin{equation}{\label{eq:3-12}}
G\Big(\frac{\omega_1}{2}\Big) - G\Big(\frac{\omega_2}{2}\Big) =
\frac{1}{8\pi} \log \left|\frac{e_1 - e_3}{e_2 - e_3}\right|.
\end{equation}

In order to compare, say, $G(\frac{1}{2}\omega_1)$ and
$G(\frac{1}{2}\omega_3)$, we may use (\ref{eq:5-5}). Let
$$
\lambda = \frac{e_3 - e_2}{e_1 - e_2}.
$$
By using $e_1 + e_2 + e_3 = 0$, we get
\begin{equation}{\label{eq:3-14}}
\frac{e_3}{e_1} = \frac{2\lambda - 1}{2 - \lambda}.
\end{equation}
It is easy to see that $|2\lambda -1| \ge |2 - \lambda|$ if and
only if $|\lambda| \ge 1$. Hence
$$
\left|\frac{e_3}{e_1}\right|\geq 1 \quad \mbox{if and only if}
\quad |\lambda|\geq 1.
$$

The same argument applies to the other two cases too and the
theorem follows.
\end{proof}

It remains to make the criterion effective in $\tau$. Recall the
modular function
$$
\lambda(\tau) = \frac{e_3 - e_2}{e_1 - e_2}.
$$
By
(\ref{eq:3-11}), we have
\begin{equation}{\label{eq:3-17}}
G\Big(\frac{\omega_3}{2}\Big) - G\Big(\frac{\omega_2}{2}\Big) =
\frac{1}{4\pi} \log |\lambda(\tau)-1|.
\end{equation}
Therefore, it is important to know when $|\lambda(\tau) - 1| = 1$.

\begin{lemma}\label{comparison}
$|\lambda(\tau) - 1| = 1$ if and only if ${\rm Re}\,\tau = \frac 1
2$.
\end{lemma}

\begin{proof} Let $\wp(z;\tau)$ be the Weierstrass $\wp$ function with
periods $1$ and $\tau$, then
$$\overline{\wp(z;\tau)}=\wp(\bar z;\bar \tau).$$
For $\tau = \frac 1 2 + ib$, $\bar\tau = 1 - \tau$ and then
$\wp(z;\bar \tau) = \wp(z;\tau)$. Thus
\begin{equation}{\label{eq:3-18}}
\wp(z;\tau) = \overline{\wp(\bar z;\tau)}\quad \mbox{for}\quad
\tau = \frac 1 2+ib.
\end{equation}
Note that if $z = \frac{1}{2}\omega_2$ then $\bar z =
\frac{1}{2}\bar\omega_2 = \frac{1}{2}(1 - \omega_2) =
\frac{1}{2}\omega_3 \pmod{\omega_1,\omega_2}$. Therefore
$$
\bar e_2 = e_3 \quad \mbox{and}\quad \bar e_1 = e_1.
$$
Since $e_1 + e_2 + e_3 = 0$, we have
\begin{equation}{\label{eq:3-19}}
{\rm Re}\, e_2 = -\frac{1}{2} e_1\quad \mbox{and}\quad {\rm Im}\,
e_2 = -{\rm Im}\, e_3.
\end{equation}
Thus
\begin{equation}{\label{eq:3-20}}
|\lambda(\tau) - 1| = \left|\frac{e_3 - e_1}{e_2 - e_1}\right| =
1.
\end{equation}

A classic result says that $\lambda'(\tau)\neq 0$ for all $\tau$.
By this and (\ref{eq:3-20}), it follows that $\lambda$ maps
$\{\,\tau \mid {\rm Re}\, \tau = \frac{1}{2}\,\}$ bijectively onto
$\{\,\lambda(\tau)\mid |\lambda(\tau) - 1| = 1\,\}$.
\end{proof}

Let $\Omega$ be the fundamental domain for $\lambda(\tau)$, i.e.,
$$
\Omega = \{\,\tau \in \mathbb{C} \mid |\tau - 1/2|\geq 1/2,\ 0\leq
{\rm Re}\, \tau\leq 1,\ {\rm Im}\,\tau > 0\,\},
$$
and let $\Omega'$ be the reflection of $\Omega$ with respect to
the imaginary axis.

Since $G(\frac{1}{2}\omega_3) < G(\frac{1}{2}\omega_2)$ for $\tau
= ib$, we conclude that for $\tau\in\Omega'\cup\Omega$,
$|\lambda(\tau) - 1| < 1$ if and only if $|{\rm Re}\, \tau| <
\frac 1 2$. Therefore for $\tau\in \Omega'\cup \Omega$,
$$
|{\rm Re}\, \tau| < \frac 1 2 \quad\mbox{if and only if}\quad
G\Big(\frac{\omega_3}{2}\Big) < G\Big(\frac{\omega_2}{2}\Big).
$$

For $|\tau|=1$, using suitable M\"obius transformations we may
obtain similar results. For example, from the definition of $\wp$,
(\ref{eq:3-18}) implies that
\begin{equation}{\label{eq:3-21}}
\bar{\wp}(z) = \Big(\frac{\tau + 1}{\bar\tau + 1}\Big)^2
\wp\Big(\frac{\tau + 1}{\bar\tau + 1}\bar z\Big)
\end{equation}
and so (compare (\ref{eq:3-1}))
\begin{equation}{\label{eq:3-22}}
G(z) = G\Big(\frac{\tau + 1}{\bar\tau + 1}\bar z\Big).
\end{equation}
Clearly, for $z = \frac 1 2 \tau$, (\ref{eq:3-22}) implies that
$G(\frac{1}{2}\omega_2) = G(\frac{1}{2}\omega_1)$. So
\begin{equation}{\label{eq:3-23}}
|\tau|=1\quad  \mbox{if and only if}\quad  \left| \frac{e_2 -
e_3}{e_1 - e_3}\right| = 1,
\end{equation}
\begin{equation}{\label{eq:3-24}}
|\tau|<1\quad \mbox{if and only if }\quad
G\Big(\frac{\omega_1}{2}\Big) < G\Big(\frac{\omega_2}{2}\Big).
\end{equation}

Similarly,
\begin{equation}{\label{eq:3-25}}
|\tau - 1| < 1\quad  \mbox{if and only if}\quad
G\Big(\frac{\omega_1}{2}\Big) < G\Big(\frac{\omega_3}{2}\Big).
\end{equation}

\section{Degeneracy Analysis of Critical Points Along
${\rm Re}\,\tau = \frac{1}{2}$} \setcounter{equation}{0}

By (\ref{eq:3-3}), the derivatives of $G$ can be computed by
\begin{equation}{\label{eq:4-1}}
\begin{split}
2\pi G_x &= {\rm Re}\, (\eta_1 t + \eta_2 s - \zeta(z)),\\
-2\pi G_y &= {\rm Im}\, (\eta_1 t + \eta_2 s - \zeta(z)).
\end{split}
\end{equation}

When $\tau = \frac 1 2 + ib$, since $\wp(z)$ is real for $z\in
\Bbb R$, $\eta_1$ is real and (\ref{eq:4-1}) becomes
\begin{equation}{\label{eq:4-2}}
\begin{split}
2\pi G_x &= \eta_1 t + \frac 1 2 \eta_1 s - {\rm Re}\,\zeta (z),\\
2\pi G_y &= {\rm Im}\,\zeta(z) + (2\pi - \eta_1 b)s.
\end{split}
\end{equation}
Thus the Hessian of $G$ is given by
\begin{equation}{\label{eq:4-3}}
\begin{split}
2\pi G_{xx} &= {\rm Re}\, \wp(z) + \eta_1,\\
2\pi G_{xy} &= - {\rm Im}\, \wp(z),\\
2\pi G_{yy} &= -\left({\rm Re}\, \wp(z) + \eta_1 -
\frac{2\pi}{b}\right).
\end{split}
\end{equation}

We first consider the point $\frac{1}{2}{\omega_1}$. The
degeneracy condition of $G$ at $\frac{1}{2}{\omega_1}$ reads as
$$
e_1 + \eta_1 = 0 \quad \mbox{or} \quad e_1 + \eta_1 -
\frac{2\pi}{b} = 0.
$$

We will use the following two inequalities (Theorem \ref{thm:1-3})
whose proofs will be given in \S8 and \S9 through theta functions:
\begin{equation}{\label{eq:4-4}}
e_1(b) + \eta_1(b)\quad \mbox{is increasing in $b$ and}
\end{equation}
\begin{equation}{\label{eq:4-5}}
e_1(b)\quad \mbox{is increasing in $b$}.
\end{equation}

\begin{lemma}\label{5-1}
There exists $b_0 < \frac{1}{2} < b_1 < \sqrt 3/2$ such that
$\frac{1}{2}{\omega_1}$ is a degenerate critical point of
$G(z;\tau)$ if and only if $b = b_0$ or $b = b_1$. Moreover,
$\frac{1}{2}{\omega_1}$ is a local minimum point of $G(z;\tau)$ if
$b\in (b_0, b_1)$ and is a saddle point of $G(z;\tau)$ if $b\in
(0,b_0)$ or $b\in (b_1,\infty)$.
\end{lemma}

\begin{proof} Let $b_0$ and $b_1$ be the zero of $e_1 + \eta_1 = 0$ and
$e_1 + \eta_1 - {2\pi}/{b} = 0$ respectively. Then Lemma \ref{5-1}
follows from the explicit expression of the Hessian of $G$ by
(\ref{eq:4-4}).
\end{proof}

Numerically we know that $b_1\approx 0.7 < {\sqrt{3}}/{2}$. Now we
analyze the behavior of $G$ near $\frac{1}{2}{\omega_1}$ for $b >
b_1$.

\begin{lemma}\label{no-real}
Suppose that $b > b_0$, then $\frac{1}{2}{\omega_1}$ is the only
critical point of $G$ along the $x$-axis.
\end{lemma}

\begin{proof} $\wp(t;\tau)$ is real if $t\in \Bbb R$. Since
$\wp'(t;\tau)\neq 0$ for $t\neq \frac{1}{2}{\omega_1}$,
$\wp'(t;\tau) < 0$ for $0 < t < \frac{1}{2}{\omega_1}$. Since $b
> b_0$, by (\ref{eq:4-3}), (\ref{eq:4-4}) and Lemma \ref{5-1},
$$
2\pi G_{xx}(t) = \wp(t) + \eta_1 > e_1 + \eta_1 > 0,
$$
which implies that $G_x(t) < G_x(\frac{1}{2}{\omega_1}) = 0$ if $0
< t < \frac{1}{2}{\omega_1}$. Hence $G$ has no critical points on
$(0, \frac{1}{2}{\omega_1})$. Since $G(z) = G(-z)$, $G$ can not
have any critical point on $(-\frac{1}{2}{\omega_1},0)$.
\end{proof}

By Lemma \ref{5-1} and the conservation of local Morse indices, we
know that $G(z;\tau)$ has two more critical points near
$\frac{1}{2}{\omega_1}$ when $b$ is close to $b_1$ and $b
> b_1$. We denote these two extra points by $z_0(\tau)$ and
$-z_0(\tau)$. In this case, $\frac{1}{2}{\omega_1}$ becomes a
saddle point and $z_0(\tau)$ and $-z_0(\tau)$ are local minimum
points. From Lemma \ref{comparison}, (\ref{eq:3-24}) and
(\ref{eq:3-25}) we know that
$$
G\Big(\frac{\omega_1}{2}\Big) < G\Big(\frac{\omega_2}{2}\Big) =
G\Big(\frac{\omega_3}{2}\Big) \quad\mbox{if}\quad b_1 < b <
{\sqrt{3}}/{2}.
$$
Thus in this region $\pm z_0(\tau)$ must exist and they turn out
to be the minimum point of $G$ since there are at most five
critical points. In fact we will see below that this is true for
all $b > b_1$ and $\frac{1}{2}\omega_2$, $\frac{1}{2}\omega_3$ are
all saddle points.

\begin{lemma}\label{5-3}
The critical point $z_0(\tau)$ is on the line ${\rm Re}\, z =
\frac{1}{2}$. Moreover, the Green function $G(z;\tau)$ is
symmetric with respect to the line ${\rm Re}\, z = \frac{1}{2}$.
\end{lemma}

\begin{proof}
Representing the torus $T$ in question by the rhombus torus with
sides $\tau$ and $\bar \tau = 1 - \tau$, then the obvious
symmetries $z \mapsto \bar z$, $z \mapsto 1 - z$ of $T$ and
Theorem \ref{uniq} show that $z_0(\tau)$ must be on the $x$-axis
or the line ${\rm Re}\, \tau = \frac{1}{2}$. The former case is
excluded by Lemma \ref{no-real}.
\end{proof}

Let $G_{\frac{1}{2}}$ be the restriction of $G$ on the line ${\rm
Re}\, z = \frac{1}{2}$, namely $G_{\frac{1}{2}}(y) = G(\frac{1}{2}
+ iy)$. Lemma \ref{5-3} implies that any critical point of
$G_{\frac{1}{2}}$ is automatically a critical point of $G$.

\begin{lemma}\label{5-4}
For $b > b_1$, $G_{\frac{1}{2}}(y)$ has exactly one critical point
in $(0,b)$. This point is necessarily a non-degenerate minimal
point.
\end{lemma}

\begin{proof} Let $z = t\omega_1 + s\omega_2$. Then ${\rm Re}\, z =
\frac{1}{2}$ is equivalent to $2t + s = 1$, which implies $\bar z
= (t + s)\omega_1 - s\omega_2 = -z$. By (\ref{eq:3-18})
$$
\wp(z;\tau) = \overline{\wp(\bar z; \tau)} = \overline{\wp
(-z;\tau)} = \overline{\wp (z;\tau)}.
$$
Hence $\wp(z;\tau)$ is real for ${\rm Re}\, z = \frac{1}{2}$. By
(\ref{eq:4-3}),
$$
2\pi G_{yy} = -\left(\wp + \eta_1 - \frac{2\pi}{b}\right).
$$
Since ${\p \wp}/{\p y} = i \wp'(z;\tau)\neq 0$ for $z\neq
\frac{1}{2}{\omega_1}$ and ${\rm Re}\, z = \frac{1}{2}$,
$G_{yy}(z)$ has at most one zero. Let $z_0(\tau)$ be the critical
point above (which exists at least for $b > b_1$ and close to
$b_1$). Then $G_y(z_0(\tau)) = G_y(\frac{1}{2}{\omega_1}) = 0$ and
so $G_{yy}(\hat z_0) = 0$ for some $\hat z_0$ in the open line
segment $(\frac{1}{2}{\omega_1}, z_0(\tau))$. Since $2\pi
G_{yy}(\frac{1}{2}{\omega_1}) = -(e_1 + \eta_1 - {2\pi}/{b}) < 0$,
we have $G_{yy}(z_0(\tau)) > 0$. Hence $z_0(\tau)$ is a
non-degenerate minimum point of $G_{\frac{1}{2}}$ as long as it
exists with $b > b_1$.

By the stability of non-degenerate minimal points (here for one
variable functions), we conclude that $z_0(\tau)$ exists for all
$b > b_1$.
\end{proof}

\begin{lemma}\label{char-rhombus-tori}
If $\wp''(z_0(\tau); \tau) = 0$ then $\tau = (1 + \sqrt{3}i)/{2}$.
\end{lemma}

\begin{proof}
Let $z_0 = t_0 \omega_1 + s_0 \omega_2$. If $\wp''(z_0;\tau) = 0$,
by the addition theorem (\ref{eq:3-6})
$$
\zeta(2z_0) = 2\zeta(z_0) = 2(t_0\eta_1 + s_0\eta_2),
$$
so $2z_0$ is also a critical point. Note that ${\rm Re}\, 2z_0 =
1$. Since $2z_0 - 1 + \omega_2$ is also a critical point with
${\rm Re}\, (2z_0 - 1 + \omega_2) = \frac{1}{2}$, we have either
$2z_0 - 1 + \omega_2 = -z_0$ or $2z_0 - 1 + \omega_2 = z_0$. The
later leads to $z_0 = 1 - \omega_2 = \omega_3$, which is not
impossible. Thus we have $2z_0 = -z_0$ in $T$ and $\wp(2z_0) =
\wp(-z_0) = \wp(z_0)$. By the addition theorem for $\wp$,
\begin{equation}{\label{addition-P}}
\wp(z_0) = \wp(2z_0) = -2\wp(z_0) + \frac 1 4 \left(
\frac{\wp''(z_0)} {\wp'(z_0)}\right)^2 = -2\wp (z_0).
\end{equation}
Therefore $\wp(z_0) = 0$. Together with $2\wp'' = 12 \wp^2 - g_2$
we find $g_2 = 0$, which is equivalent to that $\tau = (1 +
\sqrt{3}i)/{2}$.
\end{proof}

We need also the following technical lemma:

\begin{lemma}{\label{eq:4-22}}
$\wp$ maps $[\frac{1}{2}\omega_2, \frac{1}{2}\omega_3] \cup
[\frac{1}{2}(1-\omega_2), \frac{1}{2}\omega_2]$ one to one and
onto the circle $\{w\mid |w - e_1| = |e_2 - e_1|\}$, where $e_2 =
\bar e_3$, $\wp(\frac{1}{4}(\omega_2 + \omega_3)) = e_1 - |e_2 -
e_1| < 0$ and $\wp(\frac 1 4) = e_1 + |e_2 - e_1| > 0$.
\end{lemma}

Here $[\frac{1}{2}\omega_2, \frac{1}{2}\omega_3]$ means segment
$\{ \frac{1}{2}\omega_2 +t \mid 0\leq t\leq \frac 1 2\}$, and
$[\frac{1}{2}(1-\omega_2), \frac{1}{2}\omega_2]$ is $\{\frac{1}
{4} + it\mid |t|\leq \frac b 2\}$. Thus, the image of
$[\frac{1}{2}\omega_2, \frac{1}{2}\omega_3]$ is the arc connecting
$e_2$ and $e_3$ through $\wp(\frac{1}{4}(\omega_2 + \omega_3))$
and the image of $[\frac{1}{2}(1 - \omega_2),
\frac{1}{2}\omega_2]$ is the arc connecting $e_3$ and $e_2$
through $\wp(\frac 1 4)$. See Figure 4.
\begin{figure}
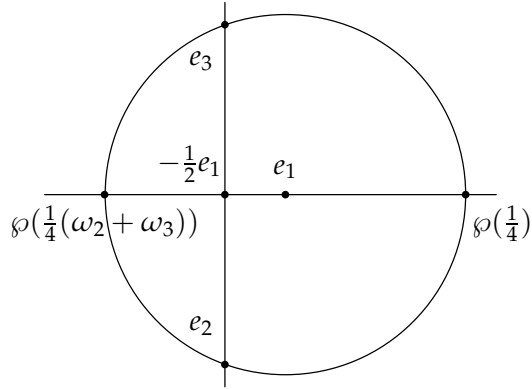

\renewcommand{\figurename}{Figure}
\begin{center}
\begin{texdraw}
  \drawdim cm  \setunitscale 0.8 \linewd 0.02
  \move(0 0) \lvec(7.5 0) \move(3 3.2) \lvec(3 -3.2)
  \move(4 0) \lcir r:3 \fcir f:0 r:0.06
  \textref h:C v:C \htext(4 0.4){$e_1$}
  \move(1 0) \fcir f:0 r:0.06
  \textref h:C v:C \htext(1 -0.5){$\wp(\frac{1}{4}
  (\omega_2 + \omega_3))$}
  \move (7 0) \fcir f:0 r:0.06
  \textref h:C v:C \htext(7.6 -0.5){$\wp(\frac{1}{4})$}
  \move (3 0) \fcir f:0 r:0.06
  \textref h:C v:C \htext(2.4 0.5){$-\frac{1}{2}e_1$}
  \move (3 2.83) \fcir f:0 r:0.06
  \textref h:C v:C \htext(2.6 2.2){$e_3$}
  \move (3 -2.83) \fcir f:0 r:0.06
  \textref h:C v:C \htext(2.6 -2.2){$e_2$}
\end{texdraw}
\end{center}
\caption{The image of the mapping $z \mapsto \wp(z)$.}
\end{figure}
Note our figure is for the case $e_1 > 0$, i.e., $b > \frac 1 2$.
In this case, the angle $\angle e_3 e_1 e_2$ is less than $\pi$.
For the case $e_1 < 0$, the angle is greater than $\pi$.

\begin{proof}
First we have
\begin{equation}{\label{eq:4-29}}
4\pi \left(G(z) - G\Big(z, \frac{\omega_1}{2}\Big)\right) = \log
\left|\wp(z) - \wp\Big(\frac{\omega_1}{2}\Big)\right| + C.
\end{equation}
Since $G(\frac{1}{2}\omega_2, \frac{1}{2}\omega_1) =
G(\frac{1}{2}\omega_2 - \frac{1}{2}\omega_1) =
G(\frac{1}{2}\omega_2)$,
\begin{equation}{\label{eq:4-30}}
4\pi \left(G(z) - G\Big(z, \frac{\omega_1}{2}\Big)\right) = \log
\left|\frac{\wp(z) - e_1} {e_2 - e_1}\right|.
\end{equation}

Let $z = \frac{1}{2}\omega_2 + t$, $t\in \Bbb R$. We have
\begin{equation*}
\begin{split}
G\Big(z, \frac{\omega_1}{2}\Big) &= G\Big(\frac{\omega_2}{2} + t -
\frac{\omega_1}{2}\Big) = G\Big(\frac{\omega_1}{2} -
\frac{\omega_2}{2} - t\Big)\\
&= G\Big(\overline{\frac{\omega_2}{2} + t}\Big) = G(z).
\end{split}
\end{equation*}
By (\ref{eq:4-30}),
\begin{equation}{\label{eq:4-31}}
\left| \frac{\wp(z) - e_1}{e_2 - e_1}\right| = 1\quad
\mbox{for}\quad z = \frac{\omega_2}{2} + t.
\end{equation}
Since $\wp(z)$ is decreasing in $y$ for $z = \frac 1 2 + iy$,
$$
\wp\Big(\frac{\omega_2 + \omega_3}{4}\Big) = \wp\Big(\frac 1 2 +
\frac{ib}{2}\Big) < \wp(1/2) = e_1.
$$
Thus $\wp(\frac{1}{4}(\omega_2 + \omega_3)) = e_1 - |e_2 - e_1|$
and the image of $[\frac{1}{2}\omega_2, \frac{1}{2}\omega_3]$ is
exactly the arc on the circle $\{\, w \mid  |w-e_1|=|e_2-e_1|\,\}$
connecting $e_2$ and $e_3$ through $\wp(\frac{1}{4}(\omega_2 +
\omega_3))$. It is one to one since $\wp'(z)\neq 0$ for $z =
\frac{1}{2}\omega_2 + t$, $t\neq 0$.

Next let $z=\frac 1 4 + it$. Then we have
\begin{equation*}
\begin{split}
G\Big(z; \frac{\omega_1}{2}\Big) &= G\Big(z -
\frac{\omega_1}{2}\Big) = G\Big(-\frac 1 4 + it\Big)
= G\Big(-\frac 1 4 - it\Big)\\
&= G\Big(\frac 1 4 + it\Big) = G(z).
\end{split}
\end{equation*}
Thus by (\ref{eq:4-30}) again,
$$
|\wp(z) - e_1| = |e_2 - e_1|\quad \mbox{for}\quad z = \frac 1 4 +
it.
$$

Since $\wp(t)$ is decreasing in $t$ for $t\in (0,\frac 1 2)$,
$\wp(\frac 1 4) > \wp(\frac 1 2)$. So $\wp(\frac 1 4) = e_1 + |e_2
- e_1|$ and the image of $[\frac{1}{2}(1-\omega_2),
\frac{1}{2}\omega_2]$ is the arc of $\{\, w\mid |w - e_1| = |e_2 -
e_1|\,\}$ connecting $e_3$ and $e_2$ through $e_1 + |e_2 - e_1|$.
\end{proof}

We are ready to prove the main results of this section:

\begin{theorem}
For $b > b_1$, $\pm z_0(\tau)$ are non-degenerate (local) minimum
points of $G$. Furthermore,
\begin{equation}{\label{eq:4-9}}
0 < {\rm Im}\, z_0(\tau) < \frac b 2.
\end{equation}
\end{theorem}

\begin{proof} We want to prove $G_{xx}(z_0(\tau);\tau) > 0$ and ${\rm
Im}\, z_0(\tau) < b/2$. By (\ref{eq:4-3}),
$$
2\pi G_{xx}(z_0(\tau);\tau) = \wp(z_0(\tau);\tau) + \eta_1.
$$

Note that $e_1e_2 + e_2e_3 + e_1e_3 = |e_2|^2 - e_1^2$ and
\begin{equation}{\label{eq:4-10}}
\wp''(z) = 2(3 \wp^2(z) + |e_2|^2 - |e_1|^2).
\end{equation}
As in (\ref{eq:3-14}), in terms of $\lambda(\tau)$,
\begin{equation}{\label{eq:4-11}}
\frac{e_2}{e_1} = \frac{\lambda + 1}{\lambda - 2}.
\end{equation}
Thus
$$
\frac{d}{d\lambda}\left(\frac{e_2}{e_1}\right) =
\frac{-3}{(\lambda-2)^2}\neq 0 \quad \mbox{for} \quad b \ne \frac
1 2.
$$
From here, we have $\displaystyle{\frac{d}{db}
\left|\frac{e_2}{e_1}\right|^2\neq 0}$ for $b \ne \frac{1}{2}$,
which implies that $\displaystyle{\frac{d}{db}
\left|\frac{e_2}{e_1}\right|^2 < 0}$ for $b > \frac{1}{2}$.
Therefore at $\tau = \frac{1}{2}(1 + \sqrt{3}i)$ (where
$\wp(z_0(\tau)) = 0$) we have
\begin{equation}{\label{eq:4-12}}
\frac{d}{db} \wp''(z_0(\tau);\tau) = 2 |e_1|^2 \frac{d}{db}
\left|\frac{e_2}{e_1}\right|^2 < 0.
\end{equation}

Since $\wp''(z_0(\tau);\tau) = 0$ at $\tau = (1 + \sqrt{3}i)/{2}$,
we have $\wp''(z_0(\tau),\tau) < 0$ for $b > \sqrt{3}/{2}$ and
sufficiently close to $\sqrt{3}/{2}$. By Lemma
\ref{char-rhombus-tori}, we then have that $\wp''(z_0(\tau);\tau)
< 0$ for all $b
> \sqrt{3}/{2}$. Thus
\begin{equation}{\label{eq:4-13}}
|\wp(z_0(\tau);\tau)|\leq \sqrt{\frac 1 3 (|e_1|^2-|e_2|^2)},
\end{equation}
and
\begin{equation}{\label{eq:4-14}}
\begin{split}
\eta_1 + \wp(z_0(\tau);\tau) &\geq  \eta_1 - \sqrt{\frac 1 3
(|e_1|^2 - |e_2|^2)}\\
&> \eta_1 - \sqrt{\frac 1 4 |e_1|^2} = \eta_1-\frac 1 2 e_1
\end{split}
\end{equation}
where $|e_2|^2 = |e_1|^2/4 + |{\rm Im}\, e_2|^2 > |e_1|^2/4$ is
used (cf.\ (\ref{eq:3-19})).

Later we will show that $\eta_1 > \frac{1}{2}e_1$ always holds
(this is part of Theorem \ref{thm:1-3} to be proved in \S9,
another direct proof will be given in (\ref{eq:4-23})). Thus the
non-degeneracy of $z_0(\tau)$ for $b
> \sqrt{3}/{2}$ follows.

For $b_1 < b < \sqrt{3}/{2}$, write $z_0(\tau) = t_0(\tau)\cdot 1
+ (1 - 2t_0(\tau))\tau$. It is clear that when $b \to b_1^+$, $t_0
\to 1/2$. We claim that $t_0(\tau)
> 1/3$ if $b_1 < b < \sqrt{3}/{2}$. For if $t_0(\tau) = 1/3$, i.e.\
$z_0(\tau) = \frac{1}{3}{\omega_3}$, then $2z_0(\tau) =
-z_0(\tau)$ is also a critical point. By the addition formula of
$\zeta$ we get $\wp''(\frac{1}{3}{\omega_3}) = 0$, and then Lemma
\ref{char-rhombus-tori} implies that $b = \sqrt{3}/2$, which is a
contradiction.

Note that by (\ref{addition-P}), it is easy to see for $\tau =
\frac{1}{2} + ib$,
\begin{equation}{\label{eq:4-16}}
12 \wp\Big(\frac{\omega_3}{3};\tau\Big) = \left(
\frac{\wp''(\frac{1}{3}{\omega_3};
\tau)}{\wp'(\frac{1}{3}{\omega_3}; \tau)}\right)^2 = -\left(
\frac{\p^2_y \wp(\frac{1}{3}\omega_3; \tau)}{\p_y
\wp(\frac{1}{3}\omega_3; \tau)} \right)^2 < 0.
\end{equation}
Hence by the monotone property of $\wp(\frac{1}{2} + iy; \tau)$ in
$y \in (0, b)$, it decreases to $-\infty$ when $y \to b$ and
\begin{equation}{\label{eq:4-15}}
\wp(z_0(\tau);\tau) > \wp\!\left(\frac{\omega_3}{3}; \tau\right).
\end{equation}

Let $f(b) := \wp(\frac{1}{3}{\omega_3}; \tau) + \frac{1}{2}
e_1(\tau)$. We have $f(\frac{1}{2}) = \wp(\frac{1}{3}{\omega_3};
\frac{1}{2}(1 + i)) < 0$ and $f(\sqrt{3}/{2}) = \frac{1}{2} e_1
(\frac{1}{2}(1 + \sqrt{3}i)) > 0$. Therefore there exists a
$\tau_0 = \frac{1}{2} + ib_0$ such that $f(b_0) = 0$. For this
$\tau_0$ and at $z = \frac{1}{3}\omega_3$ we compute
\begin{equation}
\begin{split}
\wp - e_1 &= -\frac{3}{2}e_1, \quad \wp - e_2 =
-\Big(\frac{e_1}{2}
+ e_2\Big), \quad \wp - e_3 = \frac{e_1}{2} + e_2,\\
\wp'^2 &= 4(\wp - e_1)(\wp - e_2)(\wp - e_3) =
6e_1\Big(\frac{e_1}{2} + e_2\Big)^2,\\
\wp'' &= 2\sum\nolimits_{1\leq i < j\leq 3} (\wp - e_i) (\wp -
e_j) = -2\Big(\frac{e_1}{2} + e_2\Big)^2.
\end{split}
\end{equation}

Plug in these into (\ref{eq:4-16}) we solve
\begin{equation}{\label{eq:4-19}}
\frac{e_2}{e_1} = -\frac{1}{2} \pm 3i \quad \mbox{and then} \quad
\Big|\frac{e_2}{e_1}\Big|^2 = \frac{37}{4}.
\end{equation}

Numerically at $b = b_1$, ${|e_2/e_1|^2} \approx 3.126 < 37/4$. By
the decreasing property of $|e_2/e_1|$ (c.f.\ (\ref{eq:4-12})) we
find that $\tau_0$ is unique with $b_1 > b_0$. Thus $\wp (\frac 1
3 \omega_3) + \frac 1 2 e_1 = f(b) > 0$ for $b
> b_1$. Together with (\ref{eq:4-15}),
$$
\wp(z_0(\tau),\tau) + \eta_1 > \eta_1 - \frac 1 2 e_1 > 0
$$
for $b_1 < b < {\sqrt{3}}/{2}$. The completes the proof of
$G_{xx}(z_0(\tau);\tau) > 0$.

It remains to show that ${\rm Im}\,z_0(\tau) < b/2$. This has
already been proved in the case $b_1 < b < {\sqrt{3}}/{2}$ since
$t_0(\tau) > 1/3$. For $b \ge \sqrt{3}/2$, from the continuity of
$z_0(\tau)$ in $b$ and $\wp''(z_0(\tau),\tau) \le 0$ it is enough
to show that
$$
\wp''\Big(\frac{1}{2} + \frac{1}{2}bi; \tau\Big) =
\wp''\Big(\frac{1}{4}(\omega_2 + \omega_3) \Big) > 0.
$$

Since
$$
\wp'' = \frac{1}{2}\wp'^2 \left( \frac{1}{\wp - e_1} +
\frac{1}{\wp - e_2} + \frac{1}{\wp - e_3}\right),
$$
the positivity at $\frac{1}{4}(\omega_2 + \omega_3)$ follows from
$\wp'^2 < 0$ and the negativity of the right hand side via Lemma
\ref{eq:4-22}, Figure 4.
\end{proof}

Now we discuss the non-degeneracy of $G$ at $\frac{1}{2}\omega_2$
and $\frac{1}{2}\omega_3$. The local minimum property of
$z_0(\tau)$ is in fact global by

\begin{theorem}{\label{thm:5-7}}
For $\tau=\frac 1 2 + ib$, both $\frac{1}{2}\omega_2$ and
$\frac{1}{2}\omega_3$ are non-degenerate saddle points of $G$.
\end{theorem}

\begin{proof} By (\ref{eq:4-3}), we have
\begin{equation}{\label{eq:4-20}}
\begin{split}
2\pi G_{xx}\Big(\frac{\omega_2}{2}\Big) &= {\rm
Re}\,e_2 + \eta_1 = \eta_1 - \frac 1 2 e_1,\\
2\pi G_{xy}\Big(\frac{\omega_2}{2}\Big) &= -{\rm Im}\, e_2,\\
2\pi G_{yy}\Big(\frac{\omega_2}{2}\Big) &= \frac{2\pi}{b} + \frac
1 2 e_1 - \eta_1.
\end{split}
\end{equation}
Hence the non-degeneracy of $\frac{1}{2}\omega_2$ is true if
\begin{equation}{\label{eq:4-21}}
|{\rm Im}\, e_2|^2 > \Big(\eta_1 - \frac 1 2 e_1\Big)
\Big(\frac{2\pi}{b} + \frac 1 2 e_1 - \eta_1\Big).
\end{equation}

First we claim that
\begin{equation}{\label{eq:4-23}}
\eta_1 - \frac 1 2 e_1 > 0
\end{equation}
and
\begin{equation}{\label{eq:4-24}}
\frac{2\pi}{b} + \frac 1 2 e_1 - \eta_1 > 0.
\end{equation}

To prove (\ref{eq:4-23}), we have
\begin{equation}{\label{eq:4-25}}
-\eta_1=2\int_0^{\frac 1 2} {\rm Re}\, \wp\Big(\frac{\omega_2}{2}
+ t\Big)\,dt,
\end{equation}
where $\wp(\frac{1}{2}\omega_2 + t) = \wp(\frac{1}{2}\omega_2 -
t)$ is used. By Lemma \ref{eq:4-22},
$$
{\rm Re}\, \wp\Big(\frac{\omega_2}{2} + t\Big)\leq - \frac 1 2
e_1, \quad \mbox{for all $t \in (0, \frac{1}{2})$,}
$$
hence $ -\eta_1 < -\frac 1 2 e_1$. To prove (\ref{eq:4-24}), we
have
\begin{equation*}
\begin{split}
\int_{\frac{1 - \omega_2}{2}}^{\frac{\omega_2}{2}} \wp(z)\,dz &=
\zeta\Big(\frac{1 - \omega_2}{2}\Big) -
\zeta\Big(\frac{\omega_2}{2}\Big)\\
&= \frac 1 2 (\eta_1 - 2\eta_2) = i(2\pi - b\eta_1).
\end{split}
\end{equation*}
Therefore,
$$
(2\pi - b\eta_1) = \int_{-\frac b 2}^{\frac b 2} {\rm Re}\,
\wp\Big(\frac 1 4 + it\Big)\, dt > - \frac 1 2 e_1 b
$$
and the inequality (\ref{eq:4-24}) follows.

To prove (\ref{eq:4-21}), we need two more inequalities. By
(\ref{eq:4-25}),
\begin{equation}{\label{eq:4-27}}
-\eta_1 > \wp\Big(\frac{\omega_2+\omega_3}{4}\Big) = e_1 - |e_2 -
e_1|,
\end{equation}
and by Lemma \ref{eq:4-22},
\begin{equation}{\label{eq:4-28}}
(2\pi - b\eta_1) < \wp(1/4)b = (e_1 + |e_2-e_1|)b.
\end{equation}
Thus
\begin{equation*}
\begin{split}
\Big(\eta_1 - \frac 1 2 e_1\Big)\Big(\frac{2\pi}{b} + \frac 1 2
e_1 - \eta_1\Big) & < \Big(|e_2 - e_1| - \frac 3 2 e_1\Big)
\Big(|e_2 - e_1| + \frac 3 2 e_1\Big)\\
& = |e_2 - e_1|^2 - \frac{9}{4} e_1^2 = |{\rm Im}\, e_2|^2,
\end{split}
\end{equation*}
Therefore, the non-degeneracy of $G$ at $\frac{1}{2}\omega_2$ is
proved. By (\ref{eq:4-21}), $\frac{1}{2}\omega_2$ is always a
saddle point. Representing the torus $T$ by the rhombus torus with
sides $\tau$ and $\bar \tau = 1 - \tau$, then the case for
$\frac{1}{2}\omega_3$ follows from the case for
$\frac{1}{2}\omega_2$ by symmetry reasons.
\end{proof}

\section{Green Functions Via Theta Functions}
\setcounter{equation}{0}

The purpose of \S 7 to \S 9 is to prove the two fundamental
inequalities (Theorem \ref{thm:1-3}) that have been used in
previous sections. The natural setup is based on theta functions
(we take \cite{Whittaker} as our general reference). This is easy
to explain since the moduli variable $\tau$ is explicit in theta
functions and differentiations in $\tau$ is much easier than in
the Weierstrass theory. To avoid cross references, the discussions
here are independent of previous sections.

Consider a torus $T = \mathbb{C}/\Lambda$ with $\Lambda =
(\mathbb{Z} + \mathbb{Z}\tau)$, a lattice with $\tau = a + bi$, $b
> 0$. Let $q = e^{\pi i\tau}$ with $|q| = e^{-\pi b} < 1$. The
theta function $\T_1(z; \tau)$ is the exponentially convergent
series
\begin{equation}
\begin{split}
\T_1(z; \tau) &= -i\sum_{n = -\infty}^\infty (-1)^n q^{(n +
\frac{1}{2})^2}e^{(2n + 1)\pi iz}\\
&= 2\sum_{n = 0}^\infty (-1)^n q^{(n + \frac{1}{2})^2}\sin(2n +
1)\pi z.
\end{split}
\end{equation}

For simplicity we also write it as $\T_1(z)$. It is entire with
\begin{equation}
\begin{split}
\T_1(z + 1) &= - \T_1(z),\\
\T_1(z + \tau) &= -q^{-1}e^{-2\pi iz}\T_1(z),
\end{split}
\end{equation}
which has simple zeros at the lattice points (and no others). The
following heat equation is clear from the definition
\begin{equation}
\frac{\p^2\T_1}{\p z^2} = 4\pi i\frac{\p\T_1}{\p\tau}.
\end{equation}

As usual we use $z = x + iy$. Here comes the starting point:

\begin{lemma}\label{green-theta}
Up to a constant $C(\tau)$, the Green's function $G(z, w)$ for the
Laplace operator $\triangle$ on $T$ is given by
\begin{equation}
G(z, w) = -\frac{1}{2\pi}\log |\T_1(z - w)| + \frac{1}{2b}({\rm
Im}(z - w))^2 + C(\tau).
\end{equation}
\end{lemma}

\begin{proof} Let $R(z, w)$ be the right hand side. Clearly for $z \ne w$
we have $\triangle_z R(z, w) = 1/b$ which integrated over $T$
gives 1. Near $z = w$, $R(z, w)$ has the correct behavior. So it
remains to show that $R(z, w)$ is indeed a function on $T$. From
the quasi-periodicity, $R(z + 1, w) = R(z, w)$ is obvious. Also
$$
R(z + \tau, w) - R(z, w) = -\frac{1}{2\pi}\log e^{\pi b + 2\pi y}
+ \frac{1}{2b}((y + b)^2 - y^2) = 0.
$$
These properties characterize the Green's function up to a
constant.
\end{proof}

By the translation invariance of $G$, it is enough to consider $w
= 0$. Let
$$
G(z) = G(z, 0) = -\frac{1}{2\pi}\log |\T_1(z)| + \frac{1}{2b}y^2 +
C(\tau).
$$

If we represent the torus $T$ as centered at $0$, then the
symmetry $z \mapsto -z$ shows that $G(z) = G(-z)$. By
differentiation, we get $\nabla G(z) = - \nabla G(-z)$. If $-z_0 =
z_0$ in $T$, that is $2z_0 = 0 \pmod{\Lambda}$, then we get
$\nabla G(z_0) = 0$. Hence we obtain the half periods
$\frac{1}{2}$, $\frac{1}{2}{\tau}$ and $\frac{1}{2}(1 + \tau)$ as
three obvious critical points of $G(z)$ for any $T$. By computing
$\p G/\p z = \frac{1}{2}(G_x - iG_y)$ we find
\begin{corollary}
The equation of critical points $z = x + iy$ of $G(z)$ is given by
\begin{equation}\label{critical-theta}
\frac{\p G}{\p z} \equiv \frac{-1}{4\pi}\left((\log\T_1)_z + 2\pi
i \frac{y}{b}\right) = 0.
\end{equation}
\end{corollary}

\begin{remark}
With $\zeta(z) - \eta_1 z = (\log \T_1(z))_z$ understood (cf.\
(\ref{link})), this leads to alternative simple proofs to Lemma
\ref{critical} and Corollary \ref{lem:2-3}.
\end{remark}

We compute easily
\begin{equation}
\begin{split}
G_x &= -\frac{1}{2\pi}{\rm Re}\,(\log\T_1)_z,\\
G_y &= -\frac{1}{2\pi}{\rm Re}\,(\log\T_1)_z i + \frac{y}{b} =
\frac{1}{2\pi}{\rm Im}\,(\log\T_1)_z + \frac{y}{b},\\
G_{xx} &= -\frac{1}{2\pi}{\rm Re}\,(\log\T_1)_{zz},\\
G_{xy} &= +\frac{1}{2\pi}{\rm Im}\,(\log\T_1)_{zz},\\
G_{yy} &= -\frac{1}{2\pi}{\rm Re}\,(\log\T_1)_{zz} i^2 +
\frac{1}{b} = \frac{1}{2\pi}{\rm Re}\,(\log\T_1)_{zz} +
\frac{1}{b},
\end{split}
\end{equation}
and the Hessian
\begin{equation}
\begin{split}
H &= \left|\begin{matrix} G_{xx}& G_{xy}\\ G_{yx}& G_{yy}
\end{matrix}\right| \\
&= \frac{-1}{4\pi^2} \left[({\rm Re}\,(\log\T_1)_{zz})^2 +
\frac{2\pi}{b}({\rm Re}\,(\log\T_1)_{zz}) + ({\rm Im}
\,(\log\T_1)_{zz})^2\right]\\
&= \frac{-1}{4\pi^2} \left[\Big|(\log\T_1)_{zz} +
\frac{\pi}{b}\Big|^2 - \Big(\frac{\pi}{b}\Big)^2\right].
\end{split}
\end{equation}

To analyze the critical point of $G(z)$ in general, we use the
methods of continuity to connect $\tau$ to a standard model like
the square toris, that is $\tau = i$, which under the modular
group ${\rm SL}(2,\mathbb{Z})$ is equivalent to the point $\tau =
\frac{1}{2}(1 + i)$ by $\tau \mapsto 1/(1 - \tau)$. On this
special torus, there are precisely three critical points given by
the half periods (cf.\ \cite{CLW}, Lemma 2.1).

The idea is, new critical points should be born only at certain
half period points when it degenerates under the deformation in
$\tau$. The heat equation provides a bridge between the degeneracy
condition and deformations in $\tau$. In the following, we focus
on the critical point $z = \frac{1}{2}$ and analyze its degeneracy
behavior along the half line $L$ given by $\frac{1}{2} + ib$, $b
\in \mathbb{R}$.

\section{First Inequality along the Line ${\rm Re}\,\tau =
\frac{1}{2}$} \setcounter{equation}{0}

When $\tau = \frac{1}{2} + ib \in L$, we have
$$
\T_1(z) = 2\sum_{n = 0}^\infty (-1)^n e^{\pi i/8}e^{\pi i\frac{n(n
+ 1)}{2}}e^{-\pi b(n + \frac{1}{2})^2}\sin(2n + 1)\pi z,
$$
hence the important observation that
$$
e^{-\pi i/8}\T_1(z)\in \mathbb{R}\quad\hbox{when}\quad z\in
\mathbb{R}.
$$
Similarly this holds for any derivatives of $\T_1(z)$ in $z$. In
particular,
\begin{equation}
\begin{split}
(\log\T_1)_{zz} &= \frac{{\T_1}_{zz}\T_1 -
({\T_1}_{z})^2}{\T_1^2}\\
&= 4\pi i \frac{{\T_1}_\tau}{\T_1} - (\log\T_1)_z^2 = 4\pi
(\log\T_1)_b - (\log\T_1)_z^2
\end{split}
\end{equation}
is real-valued for all $z \in \mathbb{R}$ and $\tau \in L$. Here
the heat equation and the holomorphicity of $(\log\T_1)$ have been
used.

Now we focus on the critical point $z = \frac{1}{2}$. The value $z
= \frac{1}{2}$ is fixed till the end of this section. The critical
point equation implies that $(\log\T_1)_z(\frac{1}{2}) = -2\pi i
y/b = 0$ since now $y = 0$. Thus
\begin{equation}
(\log\T_1)_{zz} = 4\pi (\log\T_1)_{b}
\end{equation}
as real functions in $b$. In this case, the point $\frac{1}{2}$ is
a degenerate critical point ($H(b) = 0$) if and only if that
\begin{equation}
4\pi (\log\T_1)_b = 0\quad\hbox{or}\quad 4\pi (\log\T_1)_b +
\frac{2\pi}{b} = 0.
\end{equation}
Notice that as functions in $b > 0$,
$$
|\T_1| = e^{-\pi i/8}\T_1\Big(\frac{1}{2}; \frac{1}{2} + ib\Big) =
2\sum_{n = 0}^\infty (-1)^{\frac{n(n + 1)}{2}}e^{-\frac{1}{4}\pi
b(2n + 1)^2}\in\mathbb{R}^+.
$$
To see this, notice that the right hand side is non-zero, real and
positive for large $b$, hence positive for all $b$. Clearly
$(\log|\T_1|)_b = (\log\T_1)_b$.

\begin{theorem}
Over the line $L$, $(\log\T_1)_{bb} = (\log|\T_1|)_{bb} < 0$.
Namely that $(\log\T_1)_b$ is decreasing from positive infinity to
$-\pi/4$. Hence that $G_{xx} = 0$ and $G_{yy} = 0$ occur exactly
once on $L$ respectively.
\end{theorem}

\begin{proof} Denote $e^{-\pi b/4}$ by $h$ and $r = h^8 = e^{-2\pi b}$.
Since $(2n + 1)^2 - 1 = 4n(n + 1)$, we get
\begin{equation}
\begin{split}
|\T_1|_b &= -2\frac{\pi}{4}\sum_{n = 0}^\infty (-1)^{\frac{n(n
+ 1)}{2}}(2n + 1)^2 h^{(2n + 1)^2} \\
&= -2h\frac{\pi}{4}\sum_{n = 0}^\infty (-1)^{\frac{n(n +
1)}{2}}(2n + 1)^2 r^{\frac{n(n + 1)}{2}}\\
|\T_1|_{bb} &= 2\frac{\pi^2}{4^2}\sum_{n = 0}^\infty
(-1)^{\frac{n(n + 1)}{2}}(2n + 1)^4 h^{(2n + 1)^2}\\
&= 2h\frac{\pi^2}{4^2}\sum_{n = 0}^\infty (-1)^{\frac{n(n +
1)}{2}}(2n + 1)^4 r^{\frac{n(n + 1)}{2}}.
\end{split}
\end{equation}

Denote the arithmetic sum $n(n + 1)/2$ by $A_n$, then
\begin{equation}
\begin{split}
&(\log|\T_1|)_{bb} = \frac{|\T_1|_{bb}|\T_1| -
|\T_1|_b^2}{|\T_1|^2}\\
&= h^2\frac{\pi^2}{4}|\T_1|^{-2}\sum_{n, m = 0}^{\infty} (-1)^{A_n
+ A_m
}((2n + 1)^4 - (2n + 1)^2 (2m + 1)^2) r^{A_n + A_m}\\
&= h^2\frac{\pi^2}{4}|\T_1|^{-2}\sum_{n > m} (-1)^{A_n + A_m
}((2n + 1)^2 - (2m + 1)^2)^2 r^{A_n + A_m}\\
&= 16h^2\pi^2|\T_1|^{-2}\sum_{n > m} (-1)^{A_n + A_m
}(A_n - A_m)^2 r^{A_n + A_m}\\
&= 16h^2\pi^2|\T_1|^{-2}(-r - 9r^3 + 4r^4 + 36r^6 - 25r^7 - 9r^9 +
100r^{10} + \cdots).
\end{split}
\end{equation}

We will prove $(\log|\T_1|)_{bb} < 0$ in two steps. First we show
by direct estimate that this is true for $b \ge \frac{1}{2}$
(indeed the argument holds for $b > 0.26$). Then we derive a
functional equation for $(\log|\T_1|)_{bb}$ which implies that the
case with $0 < b \le \frac{1}{2}$ is equivalent to the case $b \ge
\frac{1}{2}$.\smallskip

Step 1: (Direct Estimate). The point is to show that in the above
expression the sum of positive (even degree) terms is small. So
let $2k \in 2\mathbb{N}$. The number of terms with degree $2k$ is
certainly no more than $2k$, so a trivial upper bound for the
positive part is given by
\begin{equation}
A = \sum_{k =2}^\infty (2k)^3 r^{2k} = 8r^4\frac{8 - 5r^2 + 4r^4 -
r^6}{(1 - r^2)^4},
\end{equation}
where the last equality is an easy exercise in power series
calculations in Calculus. For $r \le 1/5$ we compute

\begin{equation}
\begin{split}
& (-r - 9r^3 + A)(1 - r^2)^4\\
=& -r - 5r^3 + 64r^4 + 30r^5 - 40r^6 - 50r^7 + 32r^8 + 35r^9 -
8r^{10} - 9r^{11}\\
<& -r - 5r^3 + 64r^4 + 30r^5 \\
<& -5r^3 -r\Big(1 - \frac{64}{125} - \frac{30}{625}\Big) = -5r^3 -
\frac{11}{25}r < 0.
\end{split}
\end{equation}
So $(\log|\T_1|)_{bb} < 0$ for $b = -(\log r)/2\pi > (\log 5)/2\pi
\sim 0.25615$.\smallskip

Step 2: (Functional Equation). By the Lemma to be proved below, we
have for $\hat\tau = (\tau - 1)/(2\tau -1) = \hat{a} + i\hat{b}$,
it holds that
\begin{equation}
(\log\T_1)_{\hat{b}}({1}/{2};\hat\tau) = -i(1 - 2\tau) + (1 -
2\tau)^2(\log\T_1)_b({1}/{2};\tau).
\end{equation}

When $\tau = \frac{1}{2} + ib$, we have $\hat\tau = \frac{1}{2} +
\frac{i}{4b}$. As before we may then replace $\T_1$ by $|\T_1|$.
Under $\tau\to\hat\tau$, $[\frac{1}{2},\infty)$ is mapped onto
$(0,\frac{1}{2}]$ with directions reversed. Let $f(b) =
(\log|\T_1|)_{b}(\frac{1}{2},\frac{1}{2} + ib)$. Then we get
\begin{equation}
f(1/4b) = -2b - 4b^2 f(b).
\end{equation}

Plug in $b = \frac{1}{2}$ we get that $f(\frac{1}{2}) =
-\frac{1}{2}$. So $-\frac{1}{2} = f(\frac{1}{2}) > f(b)$ for $b
> \frac{1}{2}$. Then
\begin{equation}
\begin{split}
f\Big(\frac{1}{4b}\Big) &= -2b + 2b^2 + 4b^2\Big[-\frac{1}{2} -
f(b)\Big]\\
&= 2\Big[b - \frac{1}{2}\Big]^2 - \frac{1}{2} +
4b^2\Big[-\frac{1}{2} - f(b)\Big]
\end{split}
\end{equation}
is strictly increasing in $b > \frac{1}{2}$. That is, $f(b)$ is
strictly decreasing in $b\in (0,\frac{1}{2}]$. The remaining
statements are all clear.
\end{proof}

Now we prove the functional equation. For this we need to use
Jacobi's imaginary transformation formula, which explains the
modularity for certain special theta values (cf.\ p.475 in
\cite{Whittaker}). It reads that for $\tau\tau' = -1$,
\begin{equation}
\T_1(z;\tau) = -i(i\tau')^{\frac{1}{2}}e^{\pi i \tau'
z^2}\T_1(z\tau';\tau').
\end{equation}

Recall the two generators of ${\rm SL}(2,\mathbb{Z})$ are $S\tau =
-1/\tau$ and $T\tau = \tau + 1$. Since $\T_1(z;\tau + 1) = e^{\pi
i/4}\T_1(z;\tau)$, $T$ plays no role in $(\log\T_1(z;\tau))_\tau$.

\begin{lemma}
Let $\hat\tau = ST^{-2}ST^{-1}\tau = (\tau - 1)/(2\tau -1)$. Then
\begin{equation}
(\log\T_1)_{\hat\tau}({1}/{2};\hat\tau) = -(1 - 2\tau) + (1 -
2\tau)^2(\log\T_1)_\tau({1}/{2};\tau).
\end{equation}
\end{lemma}

\begin{proof} Let $\hat\tau = S\tau_1 = -1/\tau_1$, $\tau_1 =
T^{-2}\tau_2 = \tau_2 - 2$, $\tau_2 = S\tau_3 = -1/\tau_3$ and
finally $\tau_3 = T^{-1}\tau = \tau - 1$. Notice that for
$\tau\tau' = -1$ we have $d/d\tau = \tau'^2 d/d\tau'$. Then
\begin{equation}
\begin{split}
&\frac{d}{d\hat\tau}\log\T_1({1}/{2};\hat\tau)\\
&= \tau_1^2\frac{d}{d\tau_1}\Big[\log(-i\tau_1)^{\frac{1}{2}} +
\pi i
\tau_1({1}/{2})^2 + \log\T_1(\tau_1/{2}; \tau_1)\Big]\\
&= \frac{\tau_1}{2} + \frac{\pi i \tau_1^2}{4} + \tau_1^2
\frac{d}{d\tau_1}\log\T_1(\tau_1/{2};\tau_1)\\
&= \frac{\tau_2 - 2}{2} + \frac{\pi i (\tau_2 - 2)^2}{4} + (\tau_2
- 2)^2
\frac{d}{d\tau_2}\log\T_1(\tau_2/2;\tau_2)\\
&= \frac{1}{2}\Big[\frac{-1}{\tau_3} - 2\Big] + \frac{\pi i}{4}
\Big[\frac{-1}{\tau_3} - 2\Big]^2 + \Big[\frac{-1}{\tau_3} -
2\Big]^2\times\\
&\qquad\left(\frac{\tau_3}{2} + \pi i
\tau_3^2\frac{d}{d\tau_3}\left[\tau_3(\tau_2/2)^2\right] +
\tau_3^2
\frac{d}{d\tau_3}\log\T_1(\tau_2\tau_3/2;\tau_3) \right).\\
\end{split}
\end{equation}
We plug in $\tau_2\tau_3 = -1$ and $\tau_3 = \tau - 1$. It is
clear that the second and the fourth terms are cancelled out, the
first and the third terms are combined into
\begin{equation}
\frac{1}{2}\frac{1 - 2\tau}{\tau - 1}\left[1 + \frac{1 -
2\tau}{\tau - 1}(\tau - 1)\right] = -(1 - 2\tau).
\end{equation}
This proves the result.
\end{proof}

By the previous explicit computations,
\begin{equation}
\begin{split}
G_{xx}\Big(\frac{1}{2};\frac{1}{2} + ib\Big) &=
-2(\log \T_1({1}/{2}))_b,\\
G_{xy}\Big(\frac{1}{2};\frac{1}{2} + ib\Big) &= 0,\\
G_{yy}\Big(\frac{1}{2};\frac{1}{2} + ib\Big) &=
2(\log\T_1({1}/{2}))_b + \frac{1}{b}.
\end{split}
\end{equation}
A numerical computation shows that $G_{xx}(b_0) = 0$ for $b_0 =
0.35\cdots$ and $G_{yy}(b_1) = 0$ for $b_1 = 0.71\cdots$. Hence
the sign of $(G_{xx}(b), G_{yy}(b))$ is $(-,+)$, $(+,+)$ and
$(+,-)$ for $b < b_0$, $b_0 < b < b_1$ and $b_1 < b$ respectively.
That is, $\frac{1}{2}$ is a saddle point, local minimum point and
saddle point respectively. This implies that for $b = b_1 +
\epsilon > b_1$, there are more critical points near $\frac{1}{2}$
which come from the degeneracy of $\frac{1}{2}$ at $b = b_1$ and
the conservation of local Morse index.

\section{Second Inequality along the Line ${\rm Re}\,\tau =
\frac{1}{2}$} \setcounter{equation}{0}

The analysis of extra critical points split from $\frac{1}{2}$
also relies on other though similar inequalities. Recall the three
other theta functions

\begin{equation}
\begin{split}
\T_2(z;\tau) &:= \T_1\Big(z + \frac{1}{2};\tau\Big),\\
\T_4(z;\tau) &:= \sum_{n = -\infty}^{\infty} (-1)^n q^{n^2}e^{2\pi
inz} = 1 + 2\sum_{n = 1}^\infty (-1)^n q^{n^2}\cos 2n\pi z, \\
\T_3(z;\tau) &:= \T_4\Big(z + \frac{1}{2};\tau\Big) = 1 + 2\sum_{n
= 1}^\infty q^{n^2}\cos 2n\pi z = \sum_{n = -\infty}^{\infty}
q^{n^2}e^{2\pi inz}.
\end{split}
\end{equation}
It is readily seen that $\T_1(z) = -ie^{\pi iz + \pi
i\tau/4}\T_4(z + \frac{1}{2}{\tau})$. So the four theta functions
are translates of others by half periods.

We had seen previously that $(\log|\T_1(\frac{1}{2})|)_{bb} < 0$
over the line ${{\rm Re}\,\tau} = \frac{1}{2}$. This is equivalent
to that $(\log|\T_2(0)|)_{bb} < 0$. We now discuss the case for
$\T_3(0)$ and $\T_4(0)$ where the situation is reversed(!) and it
turns out the proof is easier and purely algebraic.

\begin{theorem}
Over the line ${\rm Re}\,\tau = \frac{1}{2}$, we have $\T_4(0) =
\overline{\T_3(0)}$. Moreover, $(\log|\T_3(0)|)_{b} < 0$ and
$(\log|\T_3(0)|)_{bb} > 0$.
\end{theorem}

\begin{proof} Since $q = e^{\pi i \tau} = e^{\pi i/2}e^{-\pi b} =
ie^{-\pi b}$, $q^{n^2} = i^{n^2}e^{-\pi b n^2}$, we see that
\begin{equation}
\begin{split}
\T_3(0) &= \sum_{n\in 2\mathbb{Z}} r^{n^2} +
i\sum_{m\in2\mathbb{Z} + 1} r^{m^2},\\
\T_4(0) &= \sum_{n\in 2\mathbb{Z}} r^{n^2} -
i\sum_{m\in2\mathbb{Z} + 1} r^{m^2},
\end{split}
\end{equation}
where $r = e^{-\pi b} < 1$. Then we compute directly that
\begin{equation}
\begin{split}
|\T_3(0)|^2 &= \left(\sum\nolimits_{n\in 2\mathbb{Z}}
r^{n^2}\right)^2 + \left(\sum\nolimits_{m\in 2\mathbb{Z} + 1}
r^{m^2}\right)^2\\
&= \sum_{n, n'\in 2\mathbb{Z}} r^{n^2 + n'^2} + \sum_{m, m'\in
2\mathbb{Z} + 1} r^{m^2 + m'^2} = \sum_{k\in 2\mathbb{Z}_{\ge 0}}
p_2(k)r^k,
\end{split}
\end{equation}
where $p_2(k)$ is the number of ways to represent $k$ as the
(ordered) sum of two squares of integers. Then
\begin{equation}
(|\T_3(0)|^2)_b = -\pi \sum_{k\in 2\mathbb{Z}_{\ge 0}} p_2(k)kr^k
\end{equation}
and in particular $(\log|\T_3(0)|^2)_b < 0$.

We also have
\begin{equation}
(|\T_3(0)|^2)_{bb} = \pi^2 \sum_{k\in 2\mathbb{Z}_{\ge 0}}
p_2(k)k^2r^k.
\end{equation}
Hence that
\begin{equation}
\begin{split}
&(|\T_3(0)|^2)_{bb}|\T_3(0)|^2 - (|\T_3(0)|^2)_b^2\\
&= \pi^2\sum_{k,l\in 2\mathbb{Z}_{\ge 0}}p_2(k)p_2(l)(k^2 -
kl)r^{k + l}\\
&= \pi^2\sum_{k < l}p_2(k)p_2(l)(k - l)^2r^{k + l} > 0.
\end{split}
\end{equation}
This implies that $(\log|\T_3(0)|)_{bb} > 0$. The proof is
complete. \end{proof}

Now we relate these to Weierstrass' elliptic functions. From
$(\log\sigma(z))' = \zeta(z)$, $\sigma(z)$ is entire, odd with a
simple zero on lattice points. Moreover,
\begin{equation}
\sigma(z + \omega_i) = -e^{\eta_i(z +
\frac{1}{2}{\omega_i})}\sigma(z).
\end{equation}
This is similar to the theta function transformation law, indeed
\begin{equation}
\sigma(z) = e^{\eta_1 z^2/2}\frac{\T_1(z)}{\T_1'(0)}.
\end{equation}
Hence
\begin{equation}\label{link}
\zeta(z) - \eta_1 z = \left[\log\frac{\T_1(z)}{\T_1'(0)}\right]_z
= (\log\T_1(z))_z
\end{equation}
and
\begin{equation}
\wp(z) + \eta_1 = - (\log\T_1(z))_{zz} = - 4\pi
i(\log\T_1(z))_\tau + [(\log\T_1(z))_z]^2.
\end{equation}

For $z = \frac{1}{2}$, this simplifies to $e_1 + \eta_1 = -4\pi i
(\log\T_1(\frac{1}{2}))_\tau$. Thus our first inequality simply
says that on the line ${\rm Re}\,\tau = \frac{1}{2}$,
\begin{equation}\label{eq-1}
(e_1 + \eta_1)_b =
-4\pi\left[\log\T_1\Big(\frac{1}{2}\Big)\right]_{bb} =
-4\pi(\log\T_2(0))_{bb}>  0.
\end{equation}

From the Taylor expansion of $\sigma(z)$ and $\T_1(z)$, it is
known that
\begin{equation}
\eta_1 = -\frac{2}{3!}\frac{\T_1'''(0)}{\T_1'(0)} = -\frac{4\pi
i}{3}(\log\T_1'(0))_\tau,
\end{equation}
hence
\begin{equation}
e_1 = - 4\pi i(\log\T_2(0))_\tau + \frac{4\pi
i}{3}(\log\T_1'(0))_\tau.
\end{equation}

The Jacobi Triple Product Formula (cf.\ p.490 in \cite{Whittaker})
asserts that
$$
\T_1'(0) = \pi \T_2(0)\T_3(0)\T_4(0).
$$
So
\begin{equation}\label{eq-2}
\begin{split}
\frac{1}{2}e_1 - \eta_1 &= 2\pi i
\left[\log\frac{\T_1'(0)}{\T_2(0)} \right]_\tau \\
&= 2\pi i (\log\T_3(0)\T_4(0))_\tau = 4\pi (\log|\T_3(0)|)_b.
\end{split}
\end{equation}

Our second inequality then says that on the line ${\rm Re}\,\tau =
\frac{1}{2}$, $\frac{1}{2}e_1 - \eta_1 < 0$, $(\frac{1}{2}e_1 -
\eta_1)_b
> 0$ and $\frac{1}{2}e_1 - \eta_1$ increases to zero in $b$.
Together with the first inequality (\ref{eq-1}), we find also that
$e_1$ increases in $b$.

\end{document}